\documentclass[10pt,draftclsnofoot,onecolumn]{IEEEtran}
\ifCLASSINFOpdf
  
\else
  
\fi

\usepackage{graphicx}
\usepackage{caption}
\usepackage{subcaption}
\usepackage{amsmath} 
\usepackage{amssymb}  
\usepackage{amsthm}
\usepackage{bm}
\usepackage{lipsum}
\usepackage[linesnumbered, ruled]{algorithm2e}
\usepackage{enumerate}
\usepackage{hyperref}

\usepackage{color}

\newtheorem{theorem}{Theorem}
\newtheorem{proposition}{Proposition}
\newtheorem{lemma}{Lemma}
\newtheorem{cor}{Corollary}
\newtheorem{defn}{Definition}
\newtheorem{example}{Example}
\newtheorem{assumption}{Assumption}
\newtheorem{remark}{Remark}

\newcommand{\derivative}{D^{\tiny \mbox{S}}}
\newcommand{\newder}{D^{\tiny \mbox{NS}}}
\newcommand{\Lim}[1]{\raisebox{0.5ex}{\scalebox{0.8}{$\displaystyle \lim_{#1}\;$}}}

\begin{document}
%
\title{Approximation Algorithms for Optimization of Combinatorial Dynamical Systems}

\author{Insoon~Yang,~ Samuel A.~Burden,~ Ram~Rajagopal,~ S. Shankar~Sastry,~ and Claire J.~Tomlin
\thanks{This work was supported by the NSF CPS project ActionWebs under grant number 0931843, NSF CPS project FORCES under grant number 1239166.}
\thanks{I. Yang, S. A. Burden, S. S. Sastry and C. J. Tomlin are with the Department of Electrical Engineering and Computer Sciences, University of California, Berkeley, CA 94720, USA
        {\tt\small \{iyang, sburden, sastry, tomlin\}@eecs.berkeley.edu}}%
\thanks{R. Rajagopal is with the Department of Civil and Environmental Engineering and the Department of Electrical Engineering, Stanford University, Stanford, CA 94035, USA
        {\tt\small ramr@stanford.edu}}
}

\maketitle

\IEEEpeerreviewmaketitle

\begin{abstract}
This paper considers an optimization problem for a dynamical system whose evolution depends on a collection of binary decision variables.
We develop scalable approximation algorithms with provable suboptimality bounds to provide computationally tractable solution methods even when the dimension of the system and the number of the binary variables are large.
The proposed method employs a linear approximation of the objective function such that the approximate problem is defined over the feasible space of the binary decision variables, which is a discrete set.
To define such a linear approximation, we propose two different variation methods: one uses continuous relaxation of the discrete space and the other uses convex combinations of the vector field and running payoff.
The approximate problem is a $0$--$1$ linear program, which can be  solved by existing polynomial-time exact or approximation algorithms, and does not require the solution of the dynamical system.
Furthermore, we characterize a sufficient condition ensuring the approximate solution has a provable suboptimality bound.
We show that this condition can be interpreted as the concavity of the objective function. 
The performance and utility of the proposed algorithms are demonstrated with the ON/OFF control problems of interdependent refrigeration systems.
\end{abstract}

\section{Introduction} \label{intro}

The dynamics of critical infrastructures and their system elements---for instance, electric grid infrastructure and their electric load elements---are  
\emph{interdependent}, meaning that the state of each infrastructure or its system elements influences and is influenced by the state of the others~\cite{Rinaldi2001}. 
Such dynamic interdependencies can be classified as follows: $(i)$ infrastructure--infrastructure interdependency; $(ii)$ infrastructure--system interdependency; and $(iii)$ system--system interdependency.
All three classes of interdependencies must be addressed when making decisions that improve the performance metrics, such as efficiency, resilience and reliability, of infrastructures and their system elements.
For an example of $(i)$, consider the placement of power electronic actuators, such as high-voltage direct current links, on transmission networks. 
Such placement requires consideration of the interconnected swing dynamics of transmission grid infrastructures.
As an example of $(ii)$, it is important to consider the interdependency between the  dynamics of grid frequency and those of (aggregate) loads
when selecting the set of loads  for  frequency regulation service.
Furthermore, the ON/OFF control of a large population of electric loads whose system dynamics are coupled with each other, e.g., supermarket refrigeration systems, must take into account their system-system interdependency $(iii)$.
These decision-making problems under dynamic interdependencies 
combine the combinatorial optimization problems of
network actuator placement, load subset selection and ON/OFF control 
with the time evolution of continuous system states.
Therefore, we seek decision-making techniques that unify combinatorial optimization and dynamical systems theory.

This paper examines a fundamental problem that supports such combinatorial decision-making involving dynamical systems.
Specifically, we consider an optimization problem associated with a dynamical system whose state evolution depends on binary decision variables, which we call the \emph{combinatorial dynamical system}.
%
In our problem formulation, the binary decision variables do not change over time, unlike in the optimal control or predictive control of switched systems~\cite{Branicky1998, Xu2004, Vasudevan2013a, Bemporad1999}.
Our focus is to develop scalable methods for optimizing the binary variables associated with a dynamical system when the number of the variables is too large to enumerate all possible system `modes'
and when the dimension of the system state is large.
However, the optimization problem for combinatorial dynamical system presents a computational challenge because: $(i)$ it is a $0$--$1$ nonlinear program, which is generally NP-hard~\cite{Papadimitriou1998}; and $(ii)$ it requires the solution of a system of ordinary differential equations (ODEs).
To provide a computationally tractable solution method that can address large-scale problems, we propose scalable approximation algorithms with provable suboptimality bounds.

The key idea of the proposed methods is to linearize the objective function in the feasible space of binary decision variables.
Our first contribution is to propose a linear approximation method for  nonlinear optimization of combinatorial dynamical systems.
The approximate $0$--$1$ optimization can be efficiently solved because it is a linear $0$--$1$ program
and it does not require the solution of the dynamical system.
The proposed approximation method allows us to employ polynomial-time exact or approximation algorithms including those for problems with $l_0$-norm constraints or linear inequality constraints. 
In particular, the proposed algorithms for an $l_0$-norm constrained problem are computationally more efficient than a greedy algorithm for the same problem because our algorithms are \emph{one-shot}, i.e., do not require multiple iterations.

The proposed linear approximation approach requires the \emph{derivative} of the objective function, but this is nontrivial to construct because the function's domain is a discrete space, in general.
The second contribution of this work
is to propose two different derivative concepts. The first concept uses a natural relaxation of the discrete space, whereas  for the second concept a novel relaxation method in a function space using convex combinations of the vector fields and running payoffs is developed. We refer to the former construction as the \emph{standard derivative} because it is the same as the derivative concept in continuous space,
 and the latter as the \emph{nonstandard derivative}. 
We show the existence and the uniqueness of the nonstandard derivative, and provide an adjoint-based formula for it.
The nonstandard derivative is well-defined even when the vector field and the payoff function are undefined on interpolated values of the binary decision variables.
 Because the two derivatives are different in general, 
 we can solve two instances of the approximate problem, one with the standard derivative and another with the nonstandard derivative and then choose the better solution. 

The third contribution of this paper is to characterize conditions under which the proposed algorithms have guaranteed suboptimality bounds.
We show that the concavity of the original problem gives a sufficient condition for the suboptimality bound to hold if
 the approximation is performed using the standard derivative.
On the other hand, the same concavity condition does not suffice when the nonstandard derivative is employed in the approximation.
To resolve this difficulty, we propose a reformulated problem and show that its concavity guarantees the suboptimality bound to hold.
%
We validate the performance of the proposed approximation algorithms by solving ON/OFF control problems of commercial refrigeration systems, which consume approximately 7\% of the total commercial energy consumption in the United States~\cite{DOE2012}. 

In operations research, $0$--$1$ nonlinear optimization problems have been extensively studied over the past five decades, although the problems are not generally associated with dynamical systems.
In particular, $0$--$1$ polynomial programming, in which the objective function and the constraints are polynomials in the decision variables, has attracted great attention.
Several exact methods that can transform a $0$--$1$ polynomial program into a $0$--$1$ linear program have been developed by introducing new variables that represent the cross terms in the polynomials (e.g., \cite{Watters1967, Glover1974}).
Roof duality suggests approximation methods for $0$--$1$ polynomial programs~\cite{Hammer1984}.  
It constructs the best linear function that upperbounds the objective function (in the case of maximization) by solving a dual problem. Its size can be significantly bigger than that of the primal problem because it introduces $O(m^k)$ additional variables, where $m$ and $k$ denote the number of binary variables and the degree of polynomial, respectively.
This approach is relevant to our proposed method in the sense that both methods seek a linear function that bounds the objective function.
However, the proposed method explicitly constructs such a linear function without solving any dual problems.
Furthermore, whereas all the aforementioned methods assume that the objective function is a polynomial in the decision variables, our method does not require a polynomial representation of the objective function.
This is a considerable advantage because constructing a polynomial representation of a given function, $J:\{0,1\}^m \to \mathbb{R}$, generally requires $2^m$ calculations (e.g., via  multi-linear extension~\cite{Hammer1968}).
Even when the polynomial representations of the vector field and the objective function in the decision variables, $\alpha \in \{0,1\}^m$, are given, 
a polynomial representation of the objective function in $\alpha$ is not readily available
because the state of a dynamical system is not, in general, a polynomial in $\alpha$ with a finite degree.
For more general $0$--$1$ nonlinear programs, branch-and-bound methods (e.g., \cite{Land1960}) and penalty/smoothing methods (e.g., \cite{Murray2010}) have been suggested. However, the branch-and-bound methods cannot, in general, find a solution in polynomial time. 
The penalty and smoothing methods do not provide any performance guarantee, whereas our proposed methods guarantee suboptimality bounds.

An important class of $0$--$1$ nonlinear programs is 
the minimization or the maximization of a submodular set-function, which has the property of \emph{diminishing returns}.
Unconstrained submodular function minimization can be solved in polynomial time using a convex extension (e.g., \cite{Grotschel1981}) or a combinatorial algorithm (e.g., \cite{Schrijver2000, Iwata2001}). However, constrained submodular function minimization is NP-hard in general, and approximation algorithms with performance guarantees are available only in special cases (e.g., \cite{Goel2009, Iwata2009, Jegelka2011}).
On the other hand, our proposed method can handle a  large class of linear constraints with a guaranteed suboptimality bound.
In the case of submodular function maximization, a greedy algorithm can obtain a provably near-optimal solution~\cite{Nemhauser1978}. 
As mentioned, our proposed algorithm for $l_0$-norm constrained problems has, in general, lower computational complexity than the greedy algorithm.
We also show that the concavity conditions for our proposed suboptimality bounds to hold are not equivalent to submodularity nor does either imply the other.


The rest of this paper is organized as follows.
The problem setting for the optimization of combinatorial dynamical systems is specified in Section~\ref{setting}. 
In Section~\ref{solution}, the linear approximation approach for this problem is proposed.
To achieve the linear approximation, we propose two different concepts of the derivative of the objective function. 
Furthermore, for each linear approximation, we provide a condition under which the proposed approximate problem gives a solution with a guaranteed suboptimality bound and show that the condition can be interpreted as the concavity of the objective function or that of a reformulated objective function.
In Section~\ref{algorithm}, algorithms to solve the approximate problems with several types of linear inequality constraints are suggested.
In Section~\ref{submodular}, the proposed conditions for the suboptimality bounds to hold are compared with submodularity. 
Finally, the performance and usefulness of the proposed approximation algorithms are demonstrated with ON/OFF control problems for supermarket refrigeration systems in Section~\ref{application}.

\section{Problem Setting} \label{setting}

Consider the following dynamical system in the continuous state space $\mathcal{X} \subseteq \mathbb{R}^n$:
\begin{equation}\label{ode}
\dot{x} (t) = f(x(t), \alpha), \quad x(0) = \bold{x} \in \mathcal{X},
\end{equation}
where the vector field depends on an $m$-dimensional binary vector variable $\alpha:= \{\alpha_1, \cdots, \alpha_m\} \in \{0,1\}^m$ and $f: \mathbb{R}^n \times \mathbb{R}^m \to \mathbb{R}^n$.
We call \eqref{ode} a \emph{combinatorial dynamical system} with a binary vector variable $\alpha$.
We later view $\alpha$ as a decision variable that does not change over time in a given time interval $[0,T]$.
Let 
$x^\alpha := (x_1^\alpha, \cdots, x_n^\alpha)$
denote the solution of the ordinary differential equation \eqref{ode} given $\alpha \in \{0,1\}^m$. 
We consider the following assumptions on the vector field. 
\begin{assumption} \label{a1}
For each $\alpha \in \{0,1 \}^m$, $f(\: \cdot \: , \alpha): \mathbb{R}^n \to \mathbb{R}^n$ is twice differentiable, has a continuous second derivative and is globally Lipschitz continuous in $\mathcal{X}$. 
\end{assumption}
\begin{assumption} \label{a2}
For any $\bm{x} \in \mathcal{X}$, $f(\bm{x}, \: \cdot \:): \mathbb{R}^m \to \mathbb{R}^n$ is continuously differentiable in $[0,1]^m$.
\end{assumption}
Under Assumption~\ref{a1}, the solution of \eqref{ode} satisfies the following property~(Proposition 5.6.5 in \cite{Polak1997}): for any $\alpha \in \{0,1\}^m$,
\begin{equation}\nonumber
\|x^\alpha\|_2 := \left ( \int_0^T \| x^\alpha (t) \|^2 dt \right )^{\frac{1}{2}} < \infty.
\end{equation}
In other words, $x^\alpha : [0,T] \to \mathbb{R}^n$ is such that
$x^\alpha \in L^2([0,T]; \mathbb{R}^n)$.
Furthermore, Assumption~\ref{a1} guarantees that the system admits a unique solution, which is continuous in time, for each $\alpha \in \{0, 1\}^m$.

\subsection{Optimization of Combinatorial Dynamical Systems}

Our aim is to determine the binary vector $\alpha \in \{0,1\}^m$ that maximizes the payoff (or utility) function, $J: \mathbb{R}^m \to \mathbb{R}$, associated with the dynamical system \eqref{ode}. 
More specifically, we want to solve the following combinatorial optimization problem:
\begin{subequations} \label{opt}
\begin{align}
\max_{\alpha \in \{0,1 \}^m} \quad &J(\alpha) := \int_0^T r(x^\alpha (t),\alpha) dt + q(x^\alpha (T))\label{obj}\\
\mbox{subject to} \quad 
&\bold{A} \alpha \leq \bold{b}, \label{const}
\end{align}
\end{subequations}
where $x^\alpha$ is the solution of \eqref{ode} and
$r: \mathbb{R}^n \times \mathbb{R}^m \to \mathbb{R}$ and $q: \mathbb{R}^n  \to \mathbb{R}$ are running and terminal payoff functions, respectively.
Here, $\bold{A}$ is an $l \times m$ matrix, $\bold{b}$ is an $l$-dimensional vector and the inequality constraint \eqref{const} holds entry-wise.

This optimization problem, in general, presents a computational challenge because $(i)$ it is NP-hard; and $(ii)$ it requires the solution to the system of ODEs \eqref{ode}.
Therefore, we seek a scalable approximation method that gives a suboptimal solution with a guaranteed suboptimality bound. 
The key idea of our proposed method is to take a first-order linear approximation of the objective function \eqref{obj} with respect to the binary decision variable $\alpha$.
This linear approximation should also take into account the dependency of the state on the binary decision variable.
If the payoff function in \eqref{obj} is replaced with its linear approximation, which is linear in the decision variable, 
the approximate problem is a $0$--$1$ linear optimization.
Therefore, existing polynomial-time exact and approximation algorithms for $0$--$1$ linear programs can be employed, as shown in Section~\ref{algorithm}.
To obtain the linear approximations of the payoff function $J$,
 in the following section we formulate two different derivatives of $J$ with respect to the discrete decision variable.
Furthermore, we suggest a sufficient condition under which the approximate solution has a guaranteed suboptimality bound
in Section~\ref{guarantee}.

\section{Linear Approximation for Optimization\\ of Combinatorial Dynamical Systems} \label{solution}

Suppose for a moment that the derivative of the objective function with respect to the binary decision variable is given, and that the derivative is well-defined in $\{0,1\}^m$, which is the feasible space of the decision variable.
The derivative can be used to obtain the first-order linear approximation of the objective function, i.e., for $\alpha \in \{0,1\}^m$,
\begin{equation} \label{expand}
J(\alpha) \approx J(\bar{\alpha}) + DJ(\bar{\alpha})^\top (\alpha - \bar{\alpha}).
\end{equation}
If the objective function in \eqref{opt} is substituted with the right-hand side of \eqref{expand}, then we obtain the approximate problem:
\begin{subequations} \label{app_opt}
\begin{align}
\max_{\alpha \in \{0,1\}^m} \quad &D J (\bar{\alpha})^\top \alpha \label{app_obj} \\
\mbox{subject to} \quad &\bold{A} \alpha \leq \bold{b}.
\end{align}
\end{subequations}
This approximate problem is a 0--1 linear program, which can be solved by several polynomial-time exact or approximation algorithms (see Section~\ref{algorithm}). 
We characterize a bound on the suboptimality of the approximate solution in Section~\ref{guarantee}.

\begin{figure}[tb] 
\begin{center}
\includegraphics[width =3.3in]{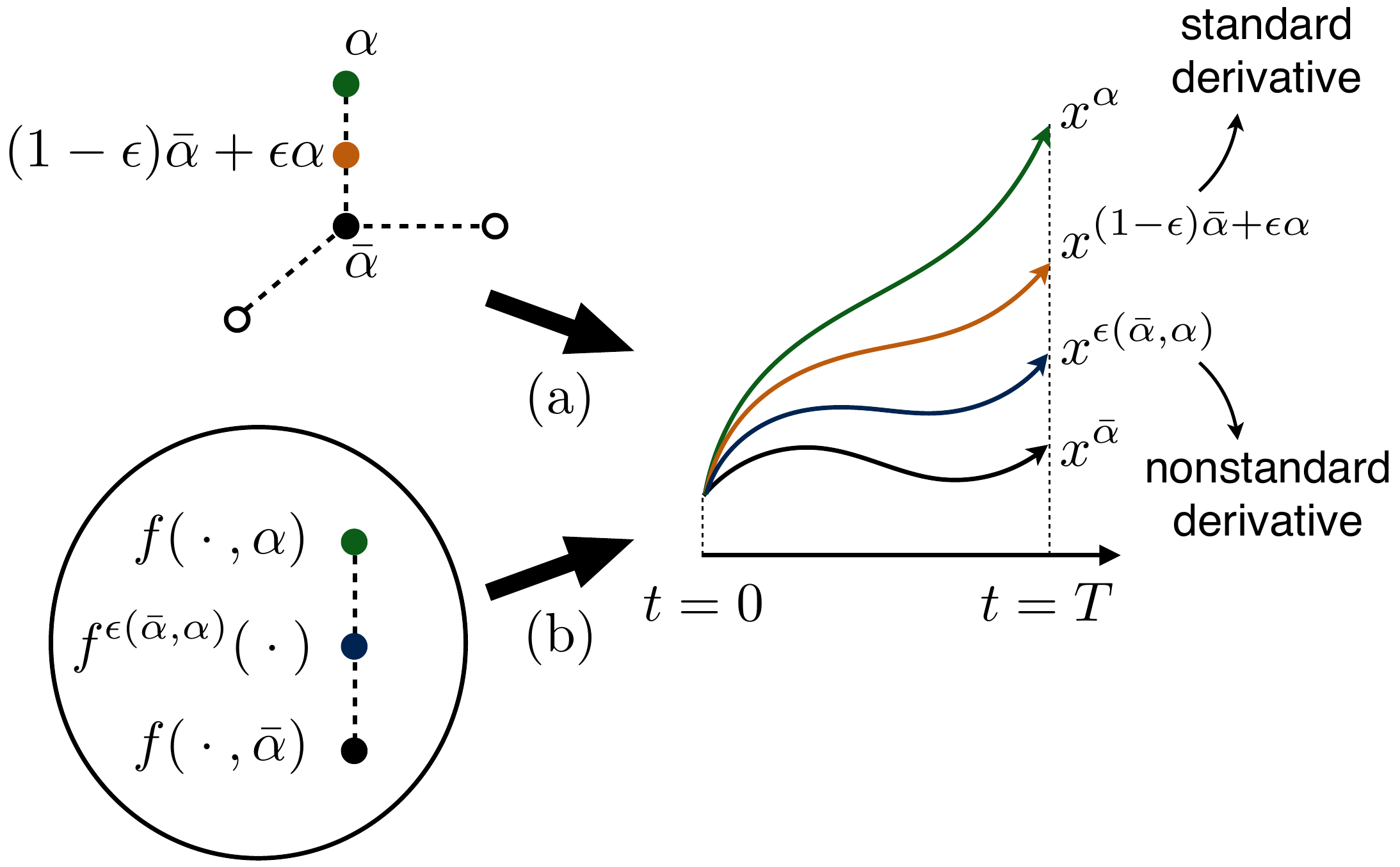}
\caption{Two variation methods: (a) the variation $(1-\epsilon) \bar{\alpha} + \epsilon \alpha$ of the binary variable produces the trajectory $x^{(1-\epsilon) \bar{\alpha} + \epsilon \alpha} (t)$, $t \in [0,T]$; and (b) the variation $f^{\epsilon(\bar{\alpha}, \alpha)}$ of the vector field, to be defined, generates another trajectory $x^{\epsilon(\bar{\alpha}, \alpha)}(t)$, $t \in [0,T]$. These two new system trajectories are used to define the standard and nonstandard derivatives, respectively.}
 \label{fig:interp}
 \end{center}
\end{figure}

We propose two different variation approaches for defining the derivatives in the discrete space $\{0,1\}^m$. The first uses the variation of  the binary decision variable in a relaxed continuous space (Fig.~\ref{fig:interp} (a)); the second uses the variation of the vector field of dynamical systems (Fig.~\ref{fig:interp} (b)). 
The first and second concepts of the derivatives are called the \emph{standard} and \emph{nonstandard} derivatives, respectively.
It is advantageous to have two different derivative concepts: we solve the approximate problem \eqref{app_opt} twice, one with the standard derivative $\derivative J$ and another with the nonstandard derivative $\newder J$ and then  choose the better solution. The one of  two approximate solutions that outperforms another is problem-dependent, in general.
We also show that the nonstandard derivative requires fewer assumptions than the standard derivative.

\subsection{Standard and Nonstandard Derivatives} \label{derivatives}


We first define the derivative of the payoff function, $J$, with respect to discrete variation of the decision variable by relaxing the discrete space $\{0, 1\}^m$ into the continuous space $\mathbb{R}^m$.
This definition of derivatives in discrete space is exactly the same as the standard definition of derivatives in continuous space.
Therefore, it requires the differentiability of the vector field and the running payoff with respect to $\alpha$.
\begin{assumption} \label{a3}
The functions $r(\: \cdot\:, \alpha): \mathbb{R}^n \to \mathbb{R}$ and $q: \mathbb{R}^n  \to \mathbb{R}$
are continuously differentiable for any $\alpha \in \{ 0, 1\}^m$.
\end{assumption}
\begin{assumption} \label{a4}
For any $\bm{x} \in \mathcal{X}$, $r(\bm{x}, \: \cdot \:): \mathbb{R}^m \to \mathbb{R}^n$ is continuously differentiable in $[0,1]^m$.
\end{assumption}
More precisely, Assumptions~\ref{a3} and \ref{a4} are needed for the standard derivative while the nonstandard derivative does not require Assumption~\ref{a4}.
Throughout this paper, we let $\bold{1}_i$ denote the $m$-dimensional vector whose $i$th entry is one and all other entries are zero.
For notational convenience, we introduce a functional, $\mathcal{J}: L^2([0,T]; \mathbb{R}^n) \times \mathbb{R}^m \to \mathbb{R}$, defined as
\begin{equation} \label{payoff_state}
\mathcal{J}(z, \beta) := \int_0^T r(z(t), \beta) dt  + q(z(T)).
\end{equation}
Note that $J(\alpha) = \mathcal{J}(x^\alpha, \alpha)$, where $x^\alpha$ is defined as the solution to the ODE \eqref{ode} with $\alpha$. 
\begin{defn} \label{defn1}
Suppose that Assumptions~\ref{a1}, \ref{a2}, \ref{a3} and \ref{a4} hold.
Given $\bar{\alpha} \in \{0,1\}^m$, the standard derivative, $\derivative J: \{0,1\}^m \to \mathbb{R}^m$, of the payoff function $J$ in \eqref{obj} 
is defined as
\begin{equation}\nonumber
[\derivative J (\bar{\alpha})]_i := \lim_{\epsilon \to 0} \frac{1}{\epsilon} \left [  
\mathcal{J} (x^{\bar{\alpha} + \epsilon \bold{1}_i}, \bar{\alpha} + \epsilon \bold{1}_i) - \mathcal{J}(x^{\bar{\alpha}}, \bar{\alpha})
 \right ]
\end{equation}
for $i=1, \cdots, m$,
where the functional $\mathcal{J}: L^2([0,T]; \mathbb{R}^n) \times \mathbb{R}^m \to \mathbb{R}$ is defined in \eqref{payoff_state} and  $x^{\bar{\alpha}}$ is the solution of \eqref{ode} with $\bar{\alpha}$.
\end{defn}
The standard derivative can be computed by direct and adjoint-based methods~\cite{Kokotovic1967, Polak1997}.
We summarize the adjoint-based method in the following proposition.
\begin{proposition} \label{prop1}
Suppose that Assumptions~\ref{a1}, \ref{a2}, \ref{a3} and \ref{a4} hold.
The derivative in Definition~\ref{defn1} can be obtained as
\begin{equation}\nonumber
\begin{split}
&\derivative J (\bar{\alpha}) = \int_0^T  \left (  \frac{\partial f(x^{\bar{\alpha}} (t), \bar{\alpha})}{\partial \bm{\alpha}}^\top \lambda^{\bar{\alpha}}(t) + \frac{\partial r (x^{\bar{\alpha}}(t), \bar{\alpha})}{\partial \bm{\alpha}}^\top \right )
dt,
\end{split}
\end{equation}
where $x^{\bar{\alpha}}$ is the solution of \eqref{ode} with $\bar{\alpha}$ and $\lambda^{\bar{\alpha}}$ solves the following adjoint system:
\begin{equation} \label{adj}
\begin{split}
- \dot{\lambda}^{\bar{\alpha}} (t) &= \frac{\partial H (x^{\bar{\alpha}}(t), \lambda^{\bar{\alpha}}(t), \bar{\alpha})}{\partial \bm{x}}^\top \\
\lambda^{\bar{\alpha}} (T) &= \frac{\partial q(x^{\bar{\alpha}}(T))}{\partial \bm{x}}^\top
\end{split}
\end{equation}
with the Hamiltonian $H: \mathbb{R}^n \times \mathbb{R}^n \times \{0,1\}^m \to \mathbb{R}$,
\begin{equation}\nonumber
H(\bm{x}, \bm{\lambda}, {\alpha}):= \bm{\lambda}^\top f(\bm{x}, {\alpha}) + r(\bm{x}, {\alpha}).
\end{equation}
\end{proposition}

We now define the derivative of the payoff function using variations in vector fields and running payoffs.
The proposed nonstandard definition of derivatives does not require
 Assumptions~\ref{a2} and \ref{a4}, i.e., the differentiability of the vector field and the running payoff with respect to $\alpha$.
Furthermore, the nonstandard derivative is well-defined even when 
 the vector field and the payoff function are not defined on the interpolated values of the binary decision variable, i.e., $f(\: \cdot \:, \alpha)$ and $r(\: \cdot \:, \alpha)$ are defined only at $\alpha \in \{0,1\}^m$.
The proposed variation procedure is as follows.
\begin{enumerate}[(i)]
\item The 0--1 vector variable $\bar{\alpha}$ in the discrete space $\{0,1\}^m$ is mapped to $x^{\bar{\alpha}}$ in the continuous metric space $L^2([0,T]; \mathbb{R}^n)$ via the original dynamical system \eqref{ode}; 
\item In $L^2([0,T]; \mathbb{R}^n)$, we construct 
a new state $x^{\epsilon (\bar{\alpha}, {\alpha})}$ as the solution to
the \emph{$\epsilon$-variational system} associated with $(\bar{\alpha}, {\alpha})$ for $\epsilon \in [0,1]$,
\begin{equation} \label{evar}
\dot{x}(t) = f^{\epsilon (\bar{\alpha}, {\alpha})} (x(t)), \quad x(0) = \bold{x} \in \mathcal{X},
\end{equation}
where the new vector field is obtained as the convex combination of the two vector fields with $\bar{\alpha}$ and $\alpha$, i.e.,
\begin{equation}\nonumber
f^{\epsilon (\bar{\alpha}, {\alpha})}(\: \cdot \:) := (1-\epsilon) f(\: \cdot \:, \bar{\alpha}) + \epsilon f(\:\cdot \:, {\alpha}).
\end{equation}
Set the distance between $\alpha$ and its $\epsilon$-variation $\epsilon (\bar{\alpha}, {\alpha})$ as $\epsilon$; and 
\item The nonstandard derivative of $J$ is defined in the following:
\end{enumerate}
\begin{defn} \label{newder_def}
Suppose that Assumptions~\ref{a1} and \ref{a3} hold. Given $\bar{\alpha} \in \{0,1\}^m$, we define the (nonstandard) derivative, $\newder J: \{0,1\}^m \to \mathbb{R}^m$ of $J$ as 
\begin{equation} \label{der2}
\begin{split}
&[\newder J (\bar{\alpha})]_i := 
\left \{ 
\begin{array}{ll}
\Lim{\epsilon  \to 0^+} \frac{1}{\epsilon}\left [  
\mathcal{J}^{\epsilon (\bar{\alpha}, \bar{\alpha} +  \bold{1}_i)} (x^{\epsilon (\bar{\alpha}, \bar{\alpha} +  \bold{1}_i)}) - \mathcal{J}(x^{\bar{\alpha}}, \bar{\alpha})
 \right ]  & \mbox{if $\bar{\alpha}_i = 0$}\\
 \Lim{\epsilon  \to 0^+} \frac{1}{\epsilon}\left [  
 \mathcal{J}(x^{\bar{\alpha}}, \bar{\alpha}) -
 \mathcal{J}^{\epsilon (\bar{\alpha}, \bar{\alpha} -  \bold{1}_i)} (x^{\epsilon (\bar{\alpha}, \bar{\alpha} -  \bold{1}_i)})
 \right ]  & \mbox{if $\bar{\alpha}_i = 1$},
\end{array}
\right.
\end{split}
\end{equation}
where $\mathcal{J}:L^2([0,T]; \mathbb{R}^n) \times \mathbb{R}^m \to \mathbb{R}$ is given by \eqref{payoff_state} and
$\mathcal{J}^{\epsilon (\bar{\alpha}, {\alpha})}: L^2([0,T]; \mathbb{R}^n) \to \mathbb{R}$ is given by
  \begin{equation} \label{varJ}
\mathcal{J}^{\epsilon (\bar{\alpha}, {\alpha})} (\: \cdot \:) :=  (1-\epsilon) \mathcal{J}(\: \cdot \:,  \bar{\alpha}) + \epsilon \mathcal{J}(\: \cdot \:, {\alpha}).
\end{equation}
Here, $x^{\bar{\alpha}}$ is the solution of \eqref{ode} with $\bar{\alpha}$ and $x^{\epsilon (\bar{\alpha}, {\alpha})}$ is the solution of \eqref{evar}.
\end{defn}
Note that we separately consider the cases with $\bar{\alpha}_i = 0$ and $\bar{\alpha}_i = 1$. This is because $\bar{\alpha} + \bold{1}_i$ is out of the feasible space of the binary decision variable when $\bar{\alpha}_i = 1$ and similarly for $\bar{\alpha} - \bold{1}_i$ when $\bar{\alpha}_i = 0$. Unlike a classical derivative with respect to continuous variable, the allowed directions for discrete variation depends on the base point $\bar{\alpha}$.
Here, the new payoff functional uses the convex combination of the running payoff because
\begin{equation} \nonumber
\mathcal{J}^{\epsilon (\bar{\alpha}, {\alpha})} (z) = 
\int_0^T (1-\epsilon) r(z, \bar{\alpha}) + \epsilon r(z, \alpha) dt + q(z(T)).
\end{equation}
The $\epsilon$-variational system is used as a continuation tool of the discrete variation from one decision variable to another.
The properties of its solution are discussed in our previous work~\cite{Yang2013} and summarized in Appendix~\ref{evar_sys}.

This nonstandard definition of derivatives raises the two following questions: $(i)$ \emph{is the nonstandard derivative well-defined?}; and $(ii)$ \emph{is there a method to compute the nonstandard derivative?} 
We answer these two questions using the adjoint system \eqref{adj} associated with the combinatorial optimization problem \eqref{opt}.

\begin{theorem} \label{dj2}
Suppose that Assumptions~\ref{a1} and \ref{a3} hold.
The nonstandard derivative $\newder J: \{0,1\}^m \to \mathbb{R}^m$ satisfies
\begin{equation} \nonumber
\begin{split}
[\newder J(\bar{\alpha})]_i :=  \int_0^T& \left ( 
f(x^{\bar{\alpha}}(t), \bar{\alpha} + \bold{1}_i) - f(x^{\bar{\alpha}}(t), \bar{\alpha}) \right )^\top \lambda^{\bar{\alpha}}(t)  + r(x^{\bar{\alpha}}(t),  \bar{\alpha}+ \bold{1}_i) - r(x^{\bar{\alpha}}(t),  \bar{\alpha} )  dt 
\end{split}
\end{equation}
if $\bar{\alpha}_i = 0$ and
\begin{equation} \nonumber
\begin{split}
[\newder J(\bar{\alpha})]_i :=  \int_0^T& \left ( 
 f(x^{\bar{\alpha}}(t), \bar{\alpha})-
 f(x^{\bar{\alpha}}(t), \bar{\alpha} - \bold{1}_i) \right )^\top \lambda^{\bar{\alpha}}(t)  +  r(x^{\bar{\alpha}}(t),  \bar{\alpha} )-
r(x^{\bar{\alpha}}(t),  \bar{\alpha}- \bold{1}_i)   dt 
\end{split}
\end{equation}
if $\bar{\alpha}_i = 1$.
Here $x^{\bar{\alpha}}$ and $\lambda^{\bar{\alpha}}$ are the solutions of \eqref{ode} and \eqref{adj} with $\bar{\alpha}$, respectively.
The derivative uniquely exists and is bounded.
\end{theorem}
The proof of Theorem~\ref{dj2} is contained in Appendix~\ref{pf_dj2}.
 The detailed comparisons  between the standard and nonstandard derivative concepts are provided in Appendix~\ref{comparison}.

\subsection{Complexity of Computing Derivatives}
\label{complexity_der}

To solve the 0--1 linear program \eqref{app_opt}, we first need to compute the standard derivative $\derivative J(\bar{\alpha})$ or the nonstandard derivative $\newder J(\bar{\alpha})$.
Recall that the dimensions of the system state and the binary decision variable are  $n$ and $m$, respectively.
Let $N_T$ be the number of time points in the time interval $[0,T]$ used to integrate the dynamical system \eqref{ode} and the adjoint system \eqref{adj}.
Then the complexity of computing the trajectories of $x^{\bar{\alpha}}$ and $\lambda^{\bar{\alpha}}$ is $O(n N_T)$ if the first-order forward Euler scheme is employed (e.g., \cite{Ascher1998}).
Note that the computation of the adjoint state trajectory $\lambda^{\bar{\alpha}}$ requires the state trajectory $x^{\bar{\alpha}}$ in $[0,T]$.
Given $x^{\bar{\alpha}}$ and $\lambda^{\bar{\alpha}}$, calculating all the entries of either the standard derivative or the nonstandard derivative requires $O( m n N_T)$ if a first-order approximation scheme for the integral over time is used. 
Therefore, the total complexity  of computing either the standard derivative or the nonstandard derivative is $O(mn N_T)$. Note that the complexity is linear in the dimension, $m$, of the decision variable $\alpha$.

\subsection{Suboptimality Bounds}\label{guarantee}

We now characterize the condition in which the solution to the approximate problem \eqref{app_opt} has a guaranteed suboptimality bound.
The suboptimality bound is obtained by showing that the optimal value of the payoff function is bounded by an affine function of the solution to the approximate problem \eqref{app_opt}.
This motivates the following concavity-like assumption:
\begin{assumption} \label{a6}
Let $\bar{\alpha} \in \{0,1\}^m$ be the point at which the original problem \eqref{opt} is linearized.
The following equality holds
\begin{equation} \label{ineq_a6}
D J(\bar{\alpha})^\top (\alpha - \bar{\alpha}) \geq J(\alpha) - J(\bar{\alpha}) \quad \forall \alpha \in\{0,1\}^m.
\end{equation}
\end{assumption}
Here, $DJ$ represents $\derivative J$ if the standard derivative used in the approximate problem, and it represents $\newder J$ if the nonstandard derivative is adopted in \eqref{app_opt}.

For notational convenience, we let $\mathcal{A}$ denote the feasible set of the optimization problem  \eqref{opt}, i.e.,
\begin{equation}\nonumber
\mathcal{A} := \{ \alpha \in \{0,1\}^m \: | \: \bold{A} \alpha \leq \bold{b}  \}.
\end{equation}
By subtracting $J(\bar{\alpha})$ from the payoff function, we normalize the payoff function such that, given $\bar{\alpha} \in \{0,1\}^m$ at which the original  problem \eqref{opt} is linearized,
\begin{equation} \nonumber
J(\bar{\alpha}) = 0.
\end{equation}
Note that $J(\alpha^{\tiny \mbox{OPT}}) \geq 0$, where $\alpha^{\tiny \mbox{OPT}}$ is a solution of the original optimization problem \eqref{opt}, if $\bar{\alpha} \in \mathcal{A}$.

\begin{theorem}[Performance Guarantee]
Suppose that Assumption~\ref{a6} holds.
Let 
\begin{equation} \label{solutions}
\begin{split} 
\alpha^{\mbox{\tiny OPT}} \in &\arg \max_{\alpha \in \mathcal{A}} \; J(\alpha),\\
\alpha^* \in &\arg \max_{\alpha \in \mathcal{A}} \; \derivative J(\bar{\alpha})^\top \alpha,\\
\hat{\alpha}^* \in &\arg \max_{\alpha \in \mathcal{A}} \; \newder J(\bar{\alpha})^\top \alpha.
\end{split}
\end{equation}
If $\derivative J(\bar{\alpha})^\top (\alpha^* - \bar{\alpha}) \neq 0$ and $\newder J(\bar{\alpha})^\top (\hat{\alpha}^* - \bar{\alpha}) \neq 0$, set
\begin{equation} \label{kappa}
\begin{split}
\rho &:= \frac{J(\alpha^*)}{\derivative J(\bar{\alpha})^\top (\alpha^* - \bar{\alpha})},\\
\hat{\rho} &:= \frac{J(\hat{\alpha}^*)}{\newder J(\bar{\alpha})^\top (\hat{\alpha}^* - \bar{\alpha})}.
\end{split}
\end{equation}
and we have the following suboptimality bounds for the solutions of the approximate problems, i.e., $\alpha^*$ and $\hat{\alpha}^*$:
\begin{equation}\label{bound1}
\begin{split}
\rho  J(\alpha^{\mbox{\tiny OPT}}) &\leq J(\alpha^*),\\
\hat{\rho}  J(\alpha^{\mbox{\tiny OPT}}) &\leq J(\hat{\alpha}^*).
\end{split}
\end{equation}
Otherwise, 
\begin{equation} \nonumber
J(\alpha^{\mbox{\tiny OPT}}) = J(\bar{\alpha}) = 0,
\end{equation}
i.e., $\bar{\alpha}$ is an optimal solution.
\end{theorem}
\begin{IEEEproof}
Due to Assumption~\ref{a6}, we have
\begin{equation} \label{i1}
J(\alpha^{\mbox{\tiny OPT}}) = J(\alpha^{\mbox{\tiny OPT}}) - J(\bar{\alpha}) \leq \derivative J(\bar{\alpha})^\top (\alpha^{\mbox{\tiny OPT}} - \bar{\alpha}).
\end{equation}
On the other hand, because ${\alpha}^* \in \arg \max_{\alpha \in \mathcal{A}} \; \derivative J(\bar{\alpha})^\top \alpha$ and $\alpha^{\mbox{\tiny OPT}} \in \mathcal{A}$, 
\begin{equation} \label{i2}
\derivative J(\bar{\alpha})^\top \alpha^{\mbox{\tiny OPT}} \leq \derivative J(\bar{\alpha})^\top {\alpha}^*.
\end{equation}
Suppose that $\derivative J(\bar{\alpha})^\top (\alpha^* - \bar{\alpha}) \neq 0$.
Combining \eqref{i1} and \eqref{i2}, we obtain the first inequality in \eqref{bound1};
the second inequality can be derived using a similar argument.
If $\derivative J(\bar{\alpha})^\top (\alpha^* - \bar{\alpha}) = 0$ or $\newder J(\bar{\alpha})^\top (\hat{\alpha}^* - \bar{\alpha}) = 0$, we have
\begin{equation} \nonumber
J(\alpha^{\mbox{\tiny OPT}}) \leq 0 = J(\bar{\alpha}).
\end{equation}
Due to the optimality of $\alpha^{\mbox{\tiny OPT}}$, the inequality must be binding. 
\end{IEEEproof}
The coefficients $\rho$ and $\hat{\rho}$ must be computed {\it a posteriori} because they require the solutions, ${\alpha}^*$ and $\hat{\alpha}^*$, respectively, of the approximate problems.
Note that $\rho$ is, in general, different from $\hat{\rho}$.
If $\bar{\alpha}$ is feasible, i.e., $\bar{\alpha} \in \mathcal{A}$, then we can improve the approximate solution by a simple post-processing that 
replaces it with $\bar{\alpha}$  if it is worse than $\bar{\alpha}$.
The payoff functions evaluated at the post-processed approximate solutions are guaranteed to be greater or than equal to zero because $J(\bar{\alpha}) = 0$.
\begin{cor}[Post-Processing] \label{bound_improved}
Suppose that Assumption~\ref{a6} holds and $\bar{\alpha} \in \mathcal{A}$. 
Let $\alpha^{\tiny \mbox{OPT}}$, $\alpha^*$ and $\hat{\alpha}^*$ be given by \eqref{solutions}.
Assume that $\derivative J(\bar{\alpha})^\top (\alpha^* - \bar{\alpha}) \neq 0$ and $\newder J(\bar{\alpha})^\top (\hat{\alpha}^* - \bar{\alpha}) \neq 0$.
Define
\begin{equation}
\begin{split}
\alpha_* &= \arg\max \{J(\alpha^*), J(\bar{\alpha})\},\\
\hat{\alpha}_* &= \arg\max \{J(\hat{\alpha}^*), J(\bar{\alpha})\}.
\end{split}
\end{equation}
and
\begin{equation} \label{sub_bounds}
\begin{split}
\rho_* &:= \max \{\rho, 0\},\\
\hat{\rho}_* &:= \max \{\hat{\rho}, 0\},
\end{split}
\end{equation}
where $\rho$ and $\hat{\rho}$ are given by \eqref{bound1}.
Then, we have the following suboptimality bounds for $\alpha_*$ and $\hat{\alpha}_*$:
\begin{equation}\label{bound2}
\begin{split}
\rho_*  J(\alpha^{\mbox{\tiny OPT}}) &\leq J(\alpha_*),\\
\hat{\rho}_*  J(\alpha^{\mbox{\tiny OPT}}) &\leq J(\hat{\alpha}_*).
\end{split}
\end{equation}
\end{cor}

The complexity of checking \eqref{ineq_a6} in Assumption~\ref{a6} for all $\alpha \in \{0,1\}^m$  increases exponentially as the dimension of the decision variable $\alpha$ increases. Therefore, we provide sufficient conditions, which are straightforward to check in some applications of interest, for Assumption~\ref{a6}.
Note that 
the inequality condition \eqref{ineq_a6} with  $DJ = \derivative J$  is equivalent to the concavity of the payoff function at $\bar{\alpha}$ if the space in which $\alpha$ lies is $[0,1]^m$ instead of $\{0, 1\}^m$. 
This observation is summarized in the following proposition.

\begin{proposition} \label{concave1}
Suppose that Assumption~\ref{a1}, \ref{a2}, \ref{a3} and \ref{a4} hold.
We also assume that the payoff function $J: \mathbb{R}^m \to \mathbb{R}$ in \eqref{obj} with $x^\alpha$ defined by \eqref{ode} is concave in $[0,1]^m$, i.e.,
\begin{equation}\nonumber
J(\alpha) := \int_0^T r(x^\alpha (t), \alpha) dt + q(x^\alpha(T)),
\end{equation}
with $x^\alpha$ satisfying
\begin{equation}\nonumber
\dot{x}^\alpha (t) = f(x^\alpha (t), \alpha), \quad x^\alpha (0)  = \bold{x} \in \mathcal{X},
\end{equation}
is concave for all $\alpha \in [0,1]^m$. 
Then, the inequality condition \eqref{ineq_a6} with $DJ = \derivative J$ holds for any $\bar{\alpha} \in \{0,1\}^m$.
\end{proposition}
Recall that we view $x^\alpha$ as a function of $\alpha$. Therefore, the concavity of $J$ is affected by how the system state depends on $\alpha$.

The inequality condition \eqref{ineq_a6} with $DJ = \newder J$ is difficult to interpret due to the nonstandard derivative.
We reformulate the dynamical system and the payoff function such that $(i)$ the standard derivative of the reformulated payoff function corresponds to the nonstandard derivative of the original payoff function and $(ii)$ the reformulated and original payoff functions have the same values at any $\alpha \in \{0,1\}^m$. Then, the concavity of the reformulated payoff function guarantees the inequality \eqref{ineq_a6}.
To be more precise, we begin by considering the following \emph{reformulated vector field and running payoff}: 
\begin{equation} \label{reform}
\begin{split}
\hat{f}(\: \cdot\: , \alpha) &:= f(\: \cdot\: , 0) + \sum_{i=1}^m \alpha_i (f (\: \cdot\: , \bold{1}_i) - f(\: \cdot\: , 0)),\\
\hat{r}(\: \cdot\: , \alpha) &:= r(\: \cdot\: , 0) + \sum_{i=1}^m \alpha_i (r (\: \cdot\: , \bold{1}_i) - r(\: \cdot\: , 0)).
\end{split}
\end{equation}
In general, $\hat{f}(\: \cdot \:, \alpha)$ (resp. $\hat{r}(\: \cdot \:, \alpha)$) and $f(\: \cdot \:, \alpha)$ (resp. $r(\: \cdot \:, \alpha)$) are different even when $\alpha$ is in the discrete space $\{0,1\}^m$. One can show that they are the same when $\alpha \in \{0,1\}^m$ if the following additivity assumption holds.
\begin{assumption} \label{a5}
The functions $f(\bm{x}, \:\cdot \:)$ and $r(\bm{x}, \: \cdot \:)$ are additive in the entries of $\alpha$ for all $\bm{x} \in \mathcal{X}$, i.e.,
\begin{equation}\nonumber
\begin{split}
f(\: \cdot \: , \alpha) &= f(\: \cdot \:, 0) + \sum_{i=1}^m (f(\: \cdot \:, \alpha_i \bold{1}_i) - f(\: \cdot \:, 0)),\\
r(\: \cdot \: , \alpha) &= r(\: \cdot \:, 0) + \sum_{i=1}^m (r(\: \cdot \:, \alpha_i \bold{1}_i) - r(\: \cdot \:, 0)).
\end{split}
\end{equation}
\end{assumption}
Note that these additivity conditions are less restrictive than the conditions that both of the functions are affine in $\alpha$ as shown in Example~\ref{ex2} in Appendix~\ref{comparison}.

This reformulation and Assumption~\ref{a5} play an essential role in interpreting the nontrivial inequality condition \eqref{ineq_a6} (with $DJ = \newder J$)
 as the concavity of
a reformulated payoff function, $\hat{J}$, defined in the next theorem.
The standard derivative of the reformulated payoff function is equivalent to the nonstandard derivative of the original payoff function under Assumption~\ref{a5}, i.e.,
\begin{equation} \nonumber
\derivative \hat{J} \equiv \newder J.
\end{equation}
Furthermore, the two payoff functions have the same values when $\alpha$ is in the discrete space $\{0,1\}^m$, i.e.,
\begin{equation} \nonumber
J|_{\{0,1\}^m} \equiv \hat{J} |_{\{0,1\}^m}.
\end{equation}
Therefore, the inequality condition \eqref{ineq_a6} with nonstandard derivative can be interpreted as the concavity of the reformulated payoff function.

\begin{theorem} \label{concave2}
Suppose that Assumptions~\ref{a1}, \ref{a3} and \ref{a5} hold.
Define the reformulated payoff function $\hat{J}: \mathbb{R}^m \to \mathbb{R}$ as 
\begin{equation}\nonumber
\hat{J}(\alpha) := \int_0^T \hat{r}(y^\alpha (t), \alpha) dt + q(y^\alpha(T)),
\end{equation}
with $y^\alpha$ satisfying
\begin{equation}\nonumber
\dot{y}^\alpha (t) = \hat{f}(y^\alpha (t), \alpha), \quad y^\alpha (0)  = \bold{x} \in \mathcal{X},
\end{equation}
where $\hat{f}$ and $\hat{r}$ are the reformulated vector field and running payoff, respectively, given in \eqref{reform}.
If the reformulated payoff function $\hat{J}$ is concave in $[0,1]^m$,
then the inequality condition \eqref{ineq_a6} with $DJ = \newder J$ holds for any $\bar{\alpha} \in \{0,1\}^m$.
\end{theorem}
\begin{IEEEproof}
Fix $\bm{x} \in \mathbb{R}^n$ and $i \in \{1, \cdots, m\}$. If $\alpha_i = 0$, then
\begin{equation}\nonumber
\alpha_i (f (\bm{x}, \bold{1}_i) - f(\bm{x}, 0)) = 0 =  f(\bm{x}, \alpha_i \bold{1}_i) - f(\bm{x}, 0).
\end{equation}
If $\alpha_i = 1$, then
\begin{equation}\nonumber
\alpha_i (f (\bm{x}, \bold{1}_i) - f(\bm{x}, 0)) = f(\bm{x}, \alpha_i \bold{1}_i) - f(\bm{x}, 0).
\end{equation}
On the other hand, due to Assumption~\ref{a5}, we have
\begin{equation}\nonumber
f(\bm{x}, \alpha) = f(\bm{x}, 0) + \sum_{i=1}^m (f(\bm{x}, \alpha_i \bold{1}_i) - f(\bm{x}, 0)).
\end{equation}
Therefore, $\hat{f}(\bm{x}, \alpha) = f(\bm{x}, \alpha)$ for any $\alpha \in \{0,1\}^m$. 
Using a similar argument, we can show that $\hat{r}(\bm{x}, \alpha) = r(\bm{x}, \alpha)$ for any $\alpha \in \{0,1\}^m$.
These imply that
\begin{equation} \label{rel1}
\hat{J}(\alpha) = J(\alpha) \quad \forall \alpha \in \{0,1\}^m.
\end{equation}
Furthermore, using the adjoint-based formula in Proposition~\ref{prop1} for the standard derivative of the reformulated payoff function $\hat{J}$, we obtain
\begin{equation} \label{format}
\begin{split}
[\derivative \hat{J}({\alpha})]_i := \int_0^T &\left ( 
f(x^{{\alpha}}(t), \bold{1}_i) - f(x^{{\alpha}}(t), 0) \right )^\top \lambda^{{\alpha}}(t)  + r(x^{{\alpha}}(t), \bold{1}_i) - r(x^{{\alpha}}(t), 0) dt.
 \end{split}
\end{equation}

On the other hand, under Assumption~\ref{a5}, the adjoint-based formula for the nonstandard derivative of the original payoff function $J$ can be rewritten as
\begin{equation} \nonumber
\begin{split}
&[\newder J({\alpha})]_i :=  \int_0^T \left ( 
f(x^{{\alpha}}(t), ({\alpha}_i +1) \bold{1}_i) - f(x^{{\alpha}}(t), {\alpha}_i \bold{1}_i) \right )^\top \lambda^{{\alpha}}(t)  + r(x^{{\alpha}}(t),  ({\alpha}_i +1) \bold{1}_i) - r(x^{{\alpha}}(t),  {\alpha}_i \bold{1}_i )  dt 
\end{split}
\end{equation}
if ${\alpha}_i = 0$ and
\begin{equation} \nonumber
\begin{split}
&[\newder J({\alpha})]_i :=  \int_0^T \left ( 
 f(x^{{\alpha}}(t), {\alpha}_i\bold{1}_i)-
 f(x^{{\alpha}}(t), ({\alpha}_i - 1)\bold{1}_i) \right )^\top \lambda^{{\alpha}}(t) +  r(x^{{\alpha}}(t),  {\alpha}_i \bold{1}_i )-
r(x^{{\alpha}}(t),  ({\alpha}-1) \bold{1}_i)   dt 
\end{split}
\end{equation}
if ${\alpha}_i = 1$.
Plugging ${\alpha}_i=0$ and ${\alpha}_i=1$ into the two formulae, respectively and comparing them with \eqref{format}, we conclude that 
\begin{equation}\label{rel2}
\derivative \hat{J}({\alpha}) =  \newder J({\alpha}) \quad \forall \alpha \in \{0,1\}^m.
\end{equation}

Suppose now that $\hat{J}$ is concave in $[0,1]^m$. Then, for any $\bar{\alpha}, \alpha \in \{0,1\}^m$,
\begin{equation} \nonumber
\derivative \hat{J}(\bar{\alpha})^\top (\alpha - \bar{\alpha}) \geq  \hat{J}(\alpha) - \hat{J}(\bar{\alpha}).
\end{equation}
Combining this inequality with \eqref{rel1} and \eqref{rel2}, we confirm that the inequality condition \eqref{ineq_a6} with $DJ = \newder J$ holds for any $\bar{\alpha} \in \{0,1\}^m$.
\end{IEEEproof}

\section{Algorithms} \label{algorithm}


We now propose approximation algorithms for the optimization of combinatorial dynamical systems \eqref{opt}
using the linear approximation proposed in the previous section. 
Formulating the approximate problem \eqref{app_opt} only requires the computation of the standard or nonstandard derivative with computational complexity $O(mnN_T)$ as suggested in Section \ref{complexity_der}, i.e., it is linear in the dimension of the decision variable. 
Because the approximate problem \eqref{app_opt} is a $0$--$1$ linear program, several polynomial time exact or approximation algorithms can be employed. 
 Another advantage of the proposed approximation is that the approximate problem no longer depends on the dynamical system. 
 Therefore, we do not need to compute the solution of the dynamical system once the derivative has been calculated.

We begin by proposing
an efficient algorithm for the $l_0$-norm constrained problem.
We then consider linear constraints \eqref{const}. 
Depending on the types of the linear constraints, several exact and approximation algorithms can be employed to solve the derivative-based approximate problem \eqref{app_opt}, which is a $0$--$1$ linear program. 

\subsection{$l_0$-Norm Constraints}\label{l0}

An important class of combinatorial optimization problems relevant to \eqref{opt} is to maximize the payoff function, given that the $l_0$-norm of the decision variable is bounded.
More specifically, instead of the original linear constraint \eqref{const}, we consider the constraint,
\begin{equation} \label{l_0_const}
\underline{K} \leq \| \alpha \|_0 \leq \overline{K},
\end{equation}
where $\underline{K}$ and $\overline{K}$ are given constants.
We consider the following first-order approximation of the combinatorial optimization problem:
\begin{equation} \label{l_0_app}
\begin{split}
\max_{\alpha \in \{0,1\}^m} \quad &DJ (\bar{\alpha})^\top \alpha \\
\mbox{subject to} \quad &\underline{K} \leq \| \alpha \|_0 \leq \overline{K}.
\end{split}
\end{equation}
A simple algorithm to solve \eqref{l_0_app} can be designed based on the ordering  of the entries of $DJ(\bar{\alpha})$, where $DJ$ is equal to either $\derivative J$ or $\newder J$.
Let $\bold{d} (\cdot)$ denote the map from $\{1, \cdots, m\}$ to $\{1, \cdots, m\}$ such that 
\begin{equation} \label{des}
[DJ(\bar{\alpha})]_{\bold{d}(i)} \geq [DJ(\bar{\alpha})]_{\bold{d}(j)}
\end{equation}
 for any $i, j \in \{1, \cdots, m\}$ such that $i \leq j$. Such a map can be constructed using a sorting algorithm with $O(m \log m)$ complexity (e.g., \cite{Sedgewick1983}).
Note that such a map may not be unique.
We  let $\alpha_{\bold{d}(i)} = 1$ for $i=1, \cdots, \underline{K}$.
We then assign $1$ on $\alpha_{\bold{d}(i)}$ if $[DJ(\bar{\alpha})]_{\bold{d}(i)} > 0$ and $\underline{K}+1 \leq i \leq \overline{K}$.
Therefore, the total computational complexity to solve the approximate problem \eqref{l_0_app} requires $O(mnN_T) + O(m \log m)$.
A more detailed algorithm to solve the $l_0$-norm constrained problem \eqref{l_0_app} is presented in Algorithm~\ref{algorithm:l_0}.

\begin{algorithm}

\textbf{Initialization:}\\

Given $\bar{\alpha}, \underline{K}, \overline{K}$; \\
$\alpha \leftarrow 0$;

\vspace{0.1in}
\textbf{Construction of $\bold{d}$:}\\
Compute $DJ(\bar{\alpha})$;\\
Sort the entries of $DJ(\bar{\alpha})$ in descending order;\\
Construct $\bold{d}: \{1, \cdots, m\} \to \{1, \cdots, m\}$ satisfying \eqref{des};

\vspace{0.1in}
\textbf{Solution of \eqref{l_0_app}:}\\
\For{$i = 1:\underline{K}$}{
$\alpha_{\bold{d}(i)} \leftarrow 1$;
}
$i \leftarrow \underline{K} + 1$;\\
\While{$[DJ(\bar{\alpha})]_{\bold{d}(i)} > 0$ and $i \leq \overline{K}$}{
$\alpha_{\bold{d}(i)} \leftarrow 1$;\\
$i \leftarrow i+1$;
}
\caption{Algorithm for the $l_0$-norm constrained problem \eqref{l_0_app}
}
\label{algorithm:l_0}
\end{algorithm}


\subsection{Totally Unimodular Matrix Constraints}

A totally unimodular (TU) matrix is defined as an integer matrix for which the determinant of every square non-singular sub-matrix is either $+1$ or $-1$.
TU matrices play an important role in integer programs because they are invertible over the integers (e.g., Chapter III.1. of~\cite{Nemhauser1988}).
Suppose that $\bold{A}$ is TU and $\bold{b}$ is integral. Let
\begin{equation}\nonumber
\bar{\bold{A}} := \begin{bmatrix}
\bold{A} \\
I_{m \times m}
\end{bmatrix}
\quad
\mbox{and } \; \;
\bar{\bold{b}} := \begin{bmatrix}
\bold{b}\\
\bold{1}
\end{bmatrix},
\end{equation}
where $\bold{1}$ is the $m$-dimensional vector whose entries are all $1$'s. The new matrix $\bar{\bold{A}}$ is also TU.
The approximate optimization problem is equivalent to the following integer linear program:
\begin{equation}\nonumber
\begin{split}
\max_{\alpha \in \mathbb{Z}^m} \quad &DJ (\bar{\alpha})^\top \alpha \\
\mbox{subject to} \quad &\bar{\bold{A}} \alpha \leq \bar{\bold{b}}.
\end{split}
\end{equation}
Because $\bar{\bold{A}}$ is TU and $\bold{b}$ is integral, 
the solution of this problem can be obtained as the solution to the linear program, whose feasible region is relaxed to $\mathbb{R}^m$, of the form
\begin{equation} \label{lp}
\begin{split}
\max_{\alpha \in \mathbb{R}^m} \quad &DJ (\bar{\alpha})^\top \alpha \\
\mbox{subject to} \quad &\bar{\bold{A}} \alpha \leq \bar{\bold{b}}.
\end{split}
\end{equation}
The proof of the exactness of this continuous relaxation can be found in~\cite{Nemhauser1988}.
The linear program \eqref{lp} can be solved by a simplex algorithm (e.g., \cite{Dantzig1998}), interior-point methods (e.g., \cite{Nesterov1994}), and several others.
Note that this approach does not require any rounding or thresholding of the solution because the solution of the relaxed problem lies in the original feasible space $\{0, 1\}^m$.

\subsection{General Linear Constraints}

Suppose that $l = 1$, i.e., $\bold{A} \in \mathbb{R}^{1 \times l}$ is a vector and $\bold{b} \in \mathbb{R}$ is a scalar and that all the entries of $\bold{A}$ and $[DJ(\bar{\alpha})]_i$ are non-negative.\footnote{The latter non-negativity assumption can easily be  relaxed by fixing $\alpha_j = 0$ for $j$ such that $[DJ(\bar{\alpha})]_j < 0$.}  
In this case,
the approximate problem \eqref{app_opt}  is a 0--1 \emph{knapsack problem}, which has been extensively studied in the past six decades.
A popular solution method is the greedy algorithm based on the linear programming (LP) relaxation proposed by Dantzig~\cite{Dantzig1957}, which replaces the feasible region $\{0,1\}^m$ with $[0,1]^m$.
A simple post-processing on the solution of the LP gives a $0.5$-approximate solution of the knapsack problem.
Such an approximate solution can be computed with complexity of $O(m) + O(m \log m)$ using a greedy algorithm (e.g., pp. 28--29 of~\cite{Martello1990}).
Other approximation algorithms have been proposed including a polynomial time approximation~\cite{Ibarra1975}.
$0$-$1$ knapsack problems with a  large number of variables can be exactly solved by branch-and-bound algorithms (e.g., \cite{Horowitz1974},  \cite{Balas1980}). 
Another classic exact method for knapsack problems is via dynamic programming (e.g., \cite{Toth1980}).
Several other algorithms and computational experiments can be found in the monograph~\cite{Martello1990} and the references therein.
If $l > 1$ and $\bold{A}_{i,j} \geq 0$ (and $[DJ(\bar{\alpha})]_j \geq 0$) for $i=1, \cdots, l$ and $j=1, \cdots, m$, then the approximate problem is called the \emph{multidimensional 0--1 knapsack problem}.  
Several exact and approximation algorithms have been developed and can be found in the review~\cite{Freville2004}, as well as among the references therein.
If no assumptions are imposed, i.e., the approximate problem \eqref{app_opt} with general linear inequality constraints is considered, then successive linear or semidefinite relaxation methods for a 0--1 polytope can provide approximation algorithms with suboptimality bounds~\cite{Lovasz1991, Sherali1990, Lasserre2002}.
\begin{remark}
Note that our proposed $0$--$1$ linear program approximation does not have any dynamical system constraints, while the original problem \eqref{opt} does.
This is advantageous because the approximate problem does not require any computational effort to solve the dynamical system once the standard or nonstandard derivative is calculated.
In other words, the complexity of any algorithm applied to the approximate problem is independent of the time horizon $[0,T]$ of the dynamical system or the number, $N_T$, of discretization points in $[0,T]$ used to approximate $DJ(\bar{\alpha})$.
\end{remark}


\section{Comparison with Submodularity} \label{submodular}

Submodularity of a set function has  attracted significant attention due to its usefulness in combinatorial optimization.  
As summarized in Section~\ref{intro}, several algorithms have been proposed for minimizing or maximizing a submodular function.
Its application includes sensor placement~\cite{Krause2008, Krause2011}, actuator placement (based on the controllability Grammian)~\cite{Summers2014},
network inference~\cite{Gomez-Rodriguez2012},
dynamic state estimation~\cite{Maillet2013}, and leader selection under link noise~\cite{Clark2014}.

Consider a set $\Omega$ with $m$ elements, 
$\Omega := \{1, \cdots, m\}$.
We define a set indicator function $\mathbb{I}: 2^\Omega \to \{0,1\}^m$ as
\begin{equation}\nonumber
[\mathbb{I}(X)]_i := \left \{
\begin{array}{ll}
0 & \mbox{if $i \notin X$}\\
1 & \mbox{if $i \in X$}.
\end{array}
\right.
\end{equation}
The set function $J(\mathbb{I}(\cdot)): 2^\Omega \to \mathbb{R}$ is said to be \emph{submodular} provided that 
for any $X \subset Y \subseteq \Omega$ and any $s \in \Omega \setminus Y$
\begin{equation}\nonumber
J(\mathbb{I}(X \cup \{s\})) - J(\mathbb{I}(X)) 
\leq
J(\mathbb{I}(Y \cup \{s\})) - J(\mathbb{I}(Y)).  
\end{equation}
If, in addition, it is monotone, i.e., for any $X \subset Y \subseteq \Omega$
\begin{equation}\nonumber
J(\mathbb{I}(X))  \leq J(\mathbb{I}(Y)), 
\end{equation}
then the problem of maximizing \eqref{obj} with $l_0$-norm constraint \eqref{l_0_const} admits a $(1 - 1/e)$  approximation algorithm~\cite{Nemhauser1978}.
Minimizing \eqref{obj} with the $l_0$-norm constraint is NP-hard, while several polynomial time algorithms can solve unconstrained submodular minimization problems as mentioned in Section~\ref{intro}.

Recall that the concavity of $J$ (resp. $\hat{J}$) guarantees the suboptimality bound to hold if the standard  (resp. nonstandard) derivative is employed (see Proposition \ref{concave1} and Theorem \ref{concave2}). 
We investigate \emph{sufficient conditions} for the concavity of $J$ and $\hat{J}$ and the submodularity of $J(\mathbb{I}(\cdot))$.
It turns out that
the concavity of $J$ or $\hat{J}$ does not imply the submodularity of $J(\mathbb{I}(\cdot))$;
furthermore,
the submodularity of $J(\mathbb{I}(\cdot))$ does not imply the concavity of $J$ or $\hat{J}$.

\subsection{Conditions for Concavity and Submodularity}

We begin by providing examples to show that concavity and  submodularity do not imply one another.
\begin{example}[Concavity does not imply submodularity]
Consider the following  vector field and running payoff:
\begin{equation}\nonumber
\begin{split}
f(\bm{x},\alpha) &= (\bm{x}_1 + \alpha_1+2, \bm{x}_2 + \alpha_2),\\
r(\bm{x},\alpha) &= -{({\bm{x}_1 - \bm{x}_2 })}^2.
\end{split}
\end{equation}
Then, we have
$r(x^\alpha (t), \alpha) = -{\left (\int_0^t e^{t-\tau} (\alpha_1 - \alpha_2+2) d\tau\right )}^2$.
The terminal payoff is set to $q \equiv 0$.
Since the following equalities hold
\begin{equation}\nonumber
\begin{split}
J(\mathbb{I}(\{2\})) - J(\mathbb{I}(\emptyset)) &= 3{\left (\int_0^t e^{t-\tau}  d\tau\right )}^2, \\
J(\mathbb{I}(\{1,2\})) - J(\mathbb{I}(\{1\})) &= 5{\left (\int_0^t 2e^{t-\tau}  d\tau\right )}^2,
 \end{split}
\end{equation}
$J(\mathbb{I}(\cdot))$ is not submodular.
On the other hand, $J = \hat{J}$ is concave in $\alpha \in [0,1]^2$.
\end{example}

\begin{example}[Submodularity does not imply concavity]
Suppose that all the assumptions in previous example hold except that the running payoff is given by
\begin{equation} \nonumber
\begin{split}
r(\bm{x},\alpha) &= {({\bm{x}_1 - \bm{x}_2 })}^2.
\end{split}
\end{equation}
In this case, $J = \hat{J}$ is not concave in $\alpha$, while $J(\mathbb{I}(\cdot))$ is submodular.
\end{example}

For comparison, 
we consider the case in which the vector field is linear in state and decision variable
 and the payoff function has a particular structure.
 In this case, the solution of the dynamical system is affine in the decision variable.

\begin{proposition} \label{prop_sub}
Suppose that $r$ is separable  as
\begin{equation}\nonumber
r(\bm{x}, \alpha) = r_1(\bm{x}) + r_2 (\alpha).
\end{equation}
Consider the vector field of the form
\begin{equation}\nonumber
f(\bm{x}, \alpha) = A\bm{x} + B \alpha,
\end{equation}
where $A$ is an $n\times n$ matrix and $B$ is an $n \times m$ matrix.
Then, \\
\begin{enumerate}
\item $J$ is concave if $r_1$,  $r_2$ and $q$ are concave;
\item $\hat{J}$ is concave if $r_1$ and $q$ are concave;
\item 
$J(\mathbb{I}(\cdot))$ is submodular  if: $r_2(\mathbb{I}(\cdot))$ is submodular; $r_1$ and $q$ are separable such that $r_1(\bm{x}) = \sum_{i=1}^m r_{1,i} (\bm{x}_i)$ and $q(\bm{x}) = \sum_{i=1}^m q_i (\bm{x}_i)$ with $r_{1,i}$ and $q_i$ concave for all $i$; and given $i$, for any $X \subset Y \subseteq \Omega$, either
\begin{equation} \nonumber
x_{i}^{\mathbb{I}(X)} (t) \leq x_{i}^{\mathbb{I}(Y)} (t)\quad \forall t \in [0,T],
\end{equation}
or
\begin{equation}\nonumber
x_{i}^{\mathbb{I}(X)} (t) \geq x_{i}^{\mathbb{I}(Y)} (t)\quad \forall t \in [0,T].
\end{equation}
\end{enumerate}
\end{proposition}
\begin{IEEEproof}
The ODE \eqref{ode} admits a unique solution,
\begin{equation}\nonumber
x^\alpha (t) = e^{At} \bold{x} + \int_0^t e^{A (t - \tau)} B \alpha d\tau.
\end{equation}
Therefore, $x^\alpha (t)$ is affine in $\alpha$ for all $t \in [0,T]$.
This implies that $r_1(x^\alpha (t))$ and $q(x^\alpha (t))$ are concave  in $\alpha$.
Furthermore, because $r_2$ is  concave  in $\alpha$, so is $J$.

In this linear system case, 
the reformulated vector field $\hat{f}$ in \eqref{reform} is equivalent to $f$
 and therefore the reformulated ODE admits the same solution, i.e., $y^\alpha \equiv x^\alpha$ for all $\alpha \in [0,1]$.
 The reformulated running payoff in \eqref{reform} is given by
\begin{equation} \nonumber
\begin{split}
\hat{r}(y^\alpha (t), \alpha) 
&=
r(y^\alpha (t), 0) + \sum_{i=1}^m \alpha_i (r (y^\alpha (t), \bold{1}_i) - r(y^\alpha (t), 0))\\
&= r(y^\alpha (t), 0) + \sum_{i=1}^m \alpha_i (r_2 (\bold{1}_i) - r_2(0)).
\end{split}
\end{equation}
Therefore, it is concave in $\alpha$ and so is $\hat{J}$.
Note that it does not require the concavity of $r_2$.

Given $i \in\{1,\cdots, n\}$, we notice that for any $X \subset Y \subseteq \Omega$ and for any $s \in \Omega \setminus Y$, either 
\begin{equation} \nonumber
x_i^{\mathbb{I}(X \cup \{s\})} - x_i^{\mathbb{I}(X )} = x_i^{\mathbb{I}(Y \cup \{s\})} - x_i^{\mathbb{I}(Y)} \geq 0, \;\; x_{i}^{\mathbb{I}(X)}  \leq x_{i}^{\mathbb{I}(Y)} 
\end{equation}
or
\begin{equation} \nonumber
x_i^{\mathbb{I}(X \cup \{s\})} - x_i^{\mathbb{I}(X )} = x_i^{\mathbb{I}(Y \cup \{s\})} - x_i^{\mathbb{I}(Y)} \leq 0, \;\; x_{i}^{\mathbb{I}(X)}  \geq x_{i}^{\mathbb{I}(Y)} 
\end{equation}
For both cases, the concavity of $r_{1,i}$ implies that
\begin{equation}\nonumber
r_{1,i}(x_i^{\mathbb{I}(X \cup \{s\})}) - r_{1,i}(x_i^{\mathbb{I}(X )}) \geq r_{1,i}(x_i^{\mathbb{I}(Y \cup \{s\})})- r_{1,i}(x_i^{\mathbb{I}(Y)}).
\end{equation}
A similar inequality holds for $q_{1,i}$.
Since $r_2$ is submodular, we also have for any $X \subset Y \subseteq \Omega$ and for any $s \in \Omega \setminus Y$
\begin{equation} \nonumber
r_2({\mathbb{I}(X \cup \{s\})}) - r_2({\mathbb{I}(X )}) \geq r_2({\mathbb{I}(Y \cup \{s\})})- r_2({\mathbb{I}(Y)}).
\end{equation}
Therefore, we obtain that for any $X \subset Y \subseteq \Omega$ and for any $s \in \Omega \setminus Y$
\begin{equation}\nonumber
\begin{split}
J(\mathbb{I}(X \cup \{s\})) - J(\mathbb{I}(X))
&= \int_0^T  r(x^{\mathbb{I}(X \cup \{s\})}, \mathbb{I}(X \cup \{s\})) - r(x^{\mathbb{I}(X)}, \mathbb{I}(X)) dt + q(x^{\mathbb{I}(X \cup \{s\})}) - q(x^{\mathbb{I}(X)})\\
&\geq
\int_0^T  r(x^{\mathbb{I}(Y \cup \{s\})}, \mathbb{I}(Y \cup \{s\})) - r(x^{\mathbb{I}(Y)}, \mathbb{I}(Y)) dt + q(x^{\mathbb{I}(Y \cup \{s\})}) - q(x^{\mathbb{I}(Y)})\\
&= J(\mathbb{I}(Y \cup \{s\})) - J(\mathbb{I}(Y)),
\end{split}
\end{equation}
which implies that $J(\mathbb{I}(\cdot))$ is submodular.
\end{IEEEproof}
Note that these are not necessary but sufficient conditions.
We observe that the concavity of $\hat{J}$ does not require the concavity of $r_2$. 
However, a sufficient condition proposed in Theorem~\ref{concave2} for Assumption~\ref{a6}, which is essential for the suboptimality bound, requires the additivity of $r_2$ (Assumption~\ref{a5}) in addition to the concavity of $\hat{J}$.
In Section~\ref{application}, the payoff function of the proposed direct load control problem satisfies all the conditions and therefore it is both concave and submodular.
In nonlinear system cases, we admit that it is nontrivial to check the concavity of $J$ or $\hat{J}$ and the submodularity of $J(\mathbb{I}(\cdot))$ unless an analytical solution of the system is available.
Further studies on characterizing the conditions for the concavity and the submodularity in the case of nonlinear systems
 will be performed in the future.

\subsection{Computational Complexity}

We now compare our proposed algorithm (Algorithm~\ref{algorithm:l_0}) for the $l_0$-norm constrained problem with the greedy algorithm for maximizing a submodular function with the same constraint,
assuming that the payoff function is submodular and satisfies Assumption~\ref{a6}.
Our algorithm is \emph{one-shot} in the sense that, after computing the derivative and ordering its entries only once, the solution is obtained.
On the other hand, the greedy algorithm chooses a locally optimal solution at each stage. In other words, this iterative greedy choice approach requires one to find an entry that maximizes the increment in the current payoff at every stage. Its complexity is $O(m^2 n N_T)$, quadratic in $m$. 
Therefore, our proposed algorithm is computationally more  efficient as the number $m$ of binary decision variables grows
 because it requires $O(m n N_T) + O(m \log m)$ calculations.

\section{Application to Direct Load Control: Interdependent Refrigeration Systems} \label{application}

The performance and usefulness of the proposed algorithms are demonstrated with applications to \emph{direct load control}, which is a demand response program in electric power systems. 
An aggregator or a load serving entity that provides a direct load control program has the authority to control its customer's loads to achieve a given objective, such as demand peak shaving, regulation services, or energy arbitrage. 
When choosing a direct load control method, a key factor is the dynamic interaction of the constituent loads.
Several scheduling methods for deferrable loads (e.g., dishwashers) with non-interdependent dynamics have been proposed that account for the variability inherent in renewable energy sources
 (e.g., \cite{Gan2013, Roozbehani2014, OBrien2013}).
 Other methods are applicable to thermostatically controlled loads with nontrivial dynamics that are decoupled (i.e., are not interdependent); see, for instance,~\cite{Malhame1985, Mathieu2013, Hao2013}.
There are fewer direct load control methods applicable to nontrivial and interdependent dynamics~\cite{Oldewurtel2013, Maasoumy2014}.

As an application, we propose a new direct load control method that can accommodate a large number of supermarkets or grocery stores containing multiple refrigerator units.
We employ a linear dynamical system model for the air temperatures near the refrigerators' evaporation units.
A refrigerant is led into an evaporator unit through an inlet valve.
Let ON (resp. OFF) control indicate that the valve is open (resp. closed).
A recent experimental study on ON/OFF control approach for supermarket refrigerator evaporators can be found in~\cite{Minetto2014}.
The control set point (ON or OFF) of one evaporator unit can affect the temperatures of other units through heat transfer.
These interdependencies 
complicate the determination of the optimal ON/OFF control.
Customer requests can yield nontrivial constraints, thereby adding further complexity to the direct load control problem.
We demonstrate that the approximation algorithms in Sections~\ref{l0} and~\ref{algorithm} are suitable for solving such direct load control problems. 
In numerical experiments, the proposed algorithms achieve near-optimal performance that surpasses that of a greedy algorithm.

Suppose that an aggregator has the authority to control refrigerator evaporator unit $i \in \{1, \cdots, m\}$ in a direct load control program.
By participating in the direct load control program, the customers who own the refrigerators can use electricity with a discounted rate.
Let \emph{air room} $i$ denote the space whose temperature is controlled by refrigerator evaporator unit $i$.
Most existing methods assume that air room $i$ is separated from air room $j$ for $i \neq j$. 
We develop a direct load control method without this assumption: our approach can handle the situation in which there is direct heat transfer between rooms $i$ and $j$ and therefore, the dynamics of  air room $i$'s temperature and the dynamics of air room $j$'s temperature are interdependent.
Let $x_i(t)$ and $\theta_i(t)$ be the temperature of room $i$ and its nearby ambient temperature at time $t$, respectively. We also let $u_i (t)$ be the ON/OFF control for refrigerator $i$ at time $t \in [0,T]$, i.e.,
\begin{equation}  \nonumber
u_i (t) := \left \{
\begin{array}{ll}
1 &\mbox{ if unit $i$ is ON at time $t$}\\
0 &\mbox{ if unit $i$ is OFF at time $t$}.
\end{array}
\right.
\end{equation}
The room temperature dynamics can be modeled by the following linear dynamical system, which is called the equivalent thermal parameter (ETP) model~\cite{Sonderegger1978}:
\begin{equation} \label{etp_original}
\dot{x}_i = -a_{ii} (x_i - \theta_i) - \sum_{j = 1}^n a_{ij} (x_i - x_j) - b_i u_i.
\end{equation}
Note that $a_{ij} = 0$ if there is no direct heat transfer between room $i$ and room $j$ (i.e., they are separated from each other). 
The ETP model can be compactly rewritten as
\begin{equation}
\dot{x} = Ax + B u + \Theta,
\end{equation}
where $A$ is an $m \times m$ connectivity matrix, whose $(i,j)$--th entry is given by
\begin{equation}\nonumber
[A]_{ij} := \left \{
\begin{array}{ll}
- \sum_{k=1}^n a_{ik} & \mbox{if $i = j$}\\
a_{ij} & \mbox{otherwise}.
\end{array}
\right.
\end{equation}
$B$ is an $m \times m$ diagonal matrix, whose $i$th diagonal entry is given by $-b_i$, and $\Theta$ is an $m$ dimensional vector, whose $i$th element is given by $a_{ii} \theta_i$.
The ambient air temperature in the supermarket is chosen as $\theta_i(t) = 19.5^\circ$C for all $t$.
If there is a path from $i$ to $j$, then the dynamics are interdependent.
\begin{figure}[tb] 
\begin{center}
\includegraphics[width =3.4in]{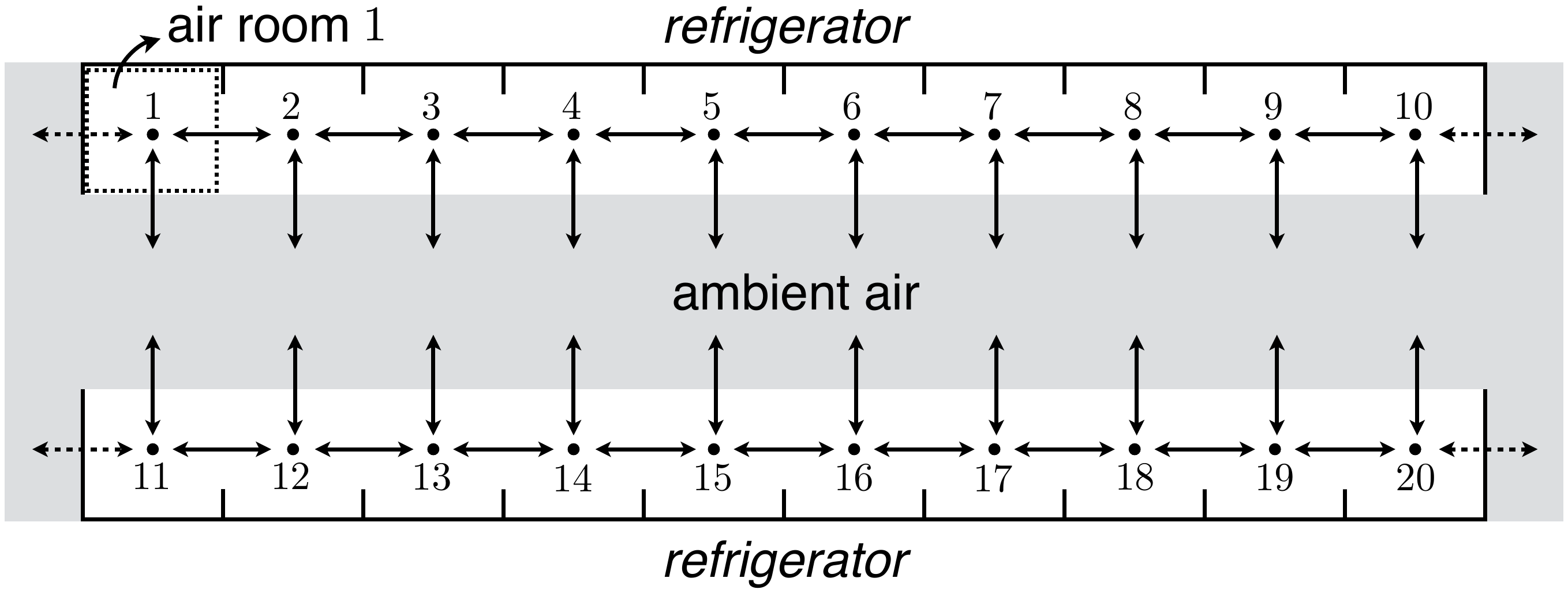}
\caption{A supermarket refrigeration system with twenty evaporator units and twenty air rooms}
 \label{fig:refrigerator}
 \end{center}
\end{figure}

Let $[\underline{\theta}_i, \overline{\theta}_i]$ be the desired temperature range for air room $i$. 
We assume that the aggregator pays the following penalty for temperature deviation to the owner of unit $i$:
\begin{equation} \nonumber
\begin{split}
&\bold{P}_i(x_i, \underline{\theta}_i, \overline{\theta}_i) := \delta_i \left [
(\underline{\theta}_i - x_i)^2 + (x_i - \overline{\theta}_i)^2 - \frac{(\underline{\theta}_i +  \overline{\theta}_i)^2}{2} 
\right ].
\end{split}
\end{equation}
The penalty is zero if $x_i = (\underline{\theta}_i + \overline{\theta}_i)/2$. 
The desired temperature range for supermarket refrigerator systems is chosen as $[\underline{\theta}_i, \overline{\theta}_i] = [0 ^\circ \mbox{C}, 4 ^\circ \mbox{C}]$. We also set $\delta_i =1$ for all $i$.

Suppose that the aggregator controls the refrigerators every $h = 15$ minutes. 
The time step can be chosen no shorter than 10 minutes because 
fast ON/OFF switching of refrigerators can cause physical failure.
The time horizon is chosen as $[10\mbox{am}, 6\mbox{pm}]$. 
The control starts from time step $1$ to time step $K=32$. Let $\alpha_i^k$ be the ON/OFF decision for unit $i$ at time step $k$. Then, the control $u_i(t)$ is set as $\alpha_i^k$ for $t \in [(k-1)h, kh)$ for $k=1, \cdots, K$.
The aggregator chooses $\alpha_i^k$ to minimize the penalty for temperature deviation 
and to provide a service to  the electric grid by following an aggregate load profile desirable to a system operator (SO).

\subsection{Case I: Target Profile} \label{target_profile}

At the beginning of time step $k$ (i.e., at time $t = (k-1)h$), the aggregator is requested by the SO to maintain the total power consumption by the $n$ refrigerator units in the target range $[\underline{y}^k, \overline{y}^k]$ (kW) for $h = 15$ minutes. In other words, the following inequality must be guaranteed:
\begin{equation} \nonumber
\underline{y}^k \leq \sum_{i=1}^n c_i \alpha_i^k \leq \overline{y}^k,
\end{equation}
where $c_i$ is the power consumption (kW) by refrigerator unit $i$ when it is in the ON state.
To reduce the energy consumption during the period of high demand,
the profile is chosen as
 $\overline{y}_k = 5500$kW for $k=9, \cdots,16$; and  $\overline{y}_k = 5000$kW otherwise in the numerical experiments.

Taking into account the penalty for temperature deviation and the constraint on the refrigerators' total power consumption,  the aggregator determines the ON/OFF control for time step $k$ as the solution of the following combinatorial optimization problem:
\begin{subequations} \label{opt_ex1}
\begin{align}
\max_{\alpha^k \in \{0,1\}^m} \: &J_k(\alpha^k):= - \int_{(k-1)h}^{kh} \sum_{i=1}^n \bold{P}_i (x_i (t), \underline{\theta}_i, \overline{\theta}_i) dt\\
\mbox{subject to} \: & \underline{y}^k \leq \sum_{i=1}^n c_i \alpha_i^k \leq \overline{y}^k \label{ineq_ex1}\\
& u(t) = \alpha^k, \quad t \in [(k-1)h, kh).
\end{align}
\end{subequations}

\begin{figure}[tb] 
\begin{center}
\includegraphics[width =3.3in]{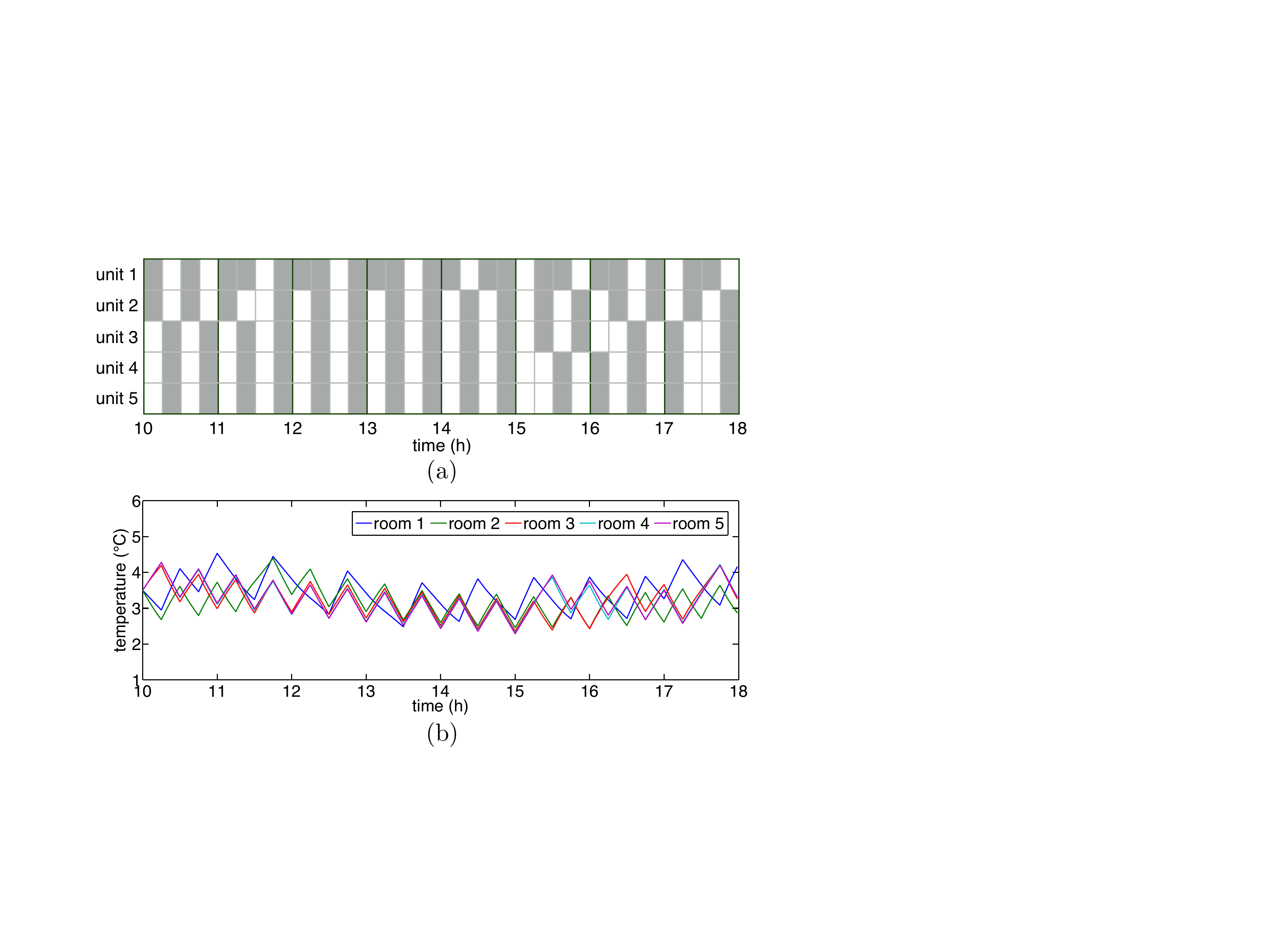}
\caption{The simulation results with approximate solution $\alpha_*^{k}$, $k=1, \cdots, 32$ for
1000 units: (a) control signals, $\alpha_{*i}^{k}$, $i=1, \cdots, 5$ (grey: ON, white: OFF) (b) controlled room temperatures, $x_i$, $i=1, \cdots, 5$.}
 \label{fig:resultTracking}
 \end{center}
\end{figure}


\begin{remark}
The inequality constraint \eqref{ineq_ex1} cannot be decomposed.
Therefore, although we can decompose the full $n$-dimensional system into $N$ subsystems such that any two subsystems are independent of each other,
the optimization problem \eqref{opt_ex1} cannot be decomposed into $N$ subproblems such that each subproblem is associated only with one of the subsystems.
\end{remark}

\begin{remark}
The payoff function, $J_k$, is concave and $J_k(\mathbb{I}(\cdot))$ is submodular due to Proposition \ref{prop_sub}. Furthermore, the standard and the nonstandard derivatives are the same in this problem due to Proposition \ref{prop_equiv}.
\end{remark}

\begin{figure}[tb] 
\begin{center}
\includegraphics[width =3.3in]{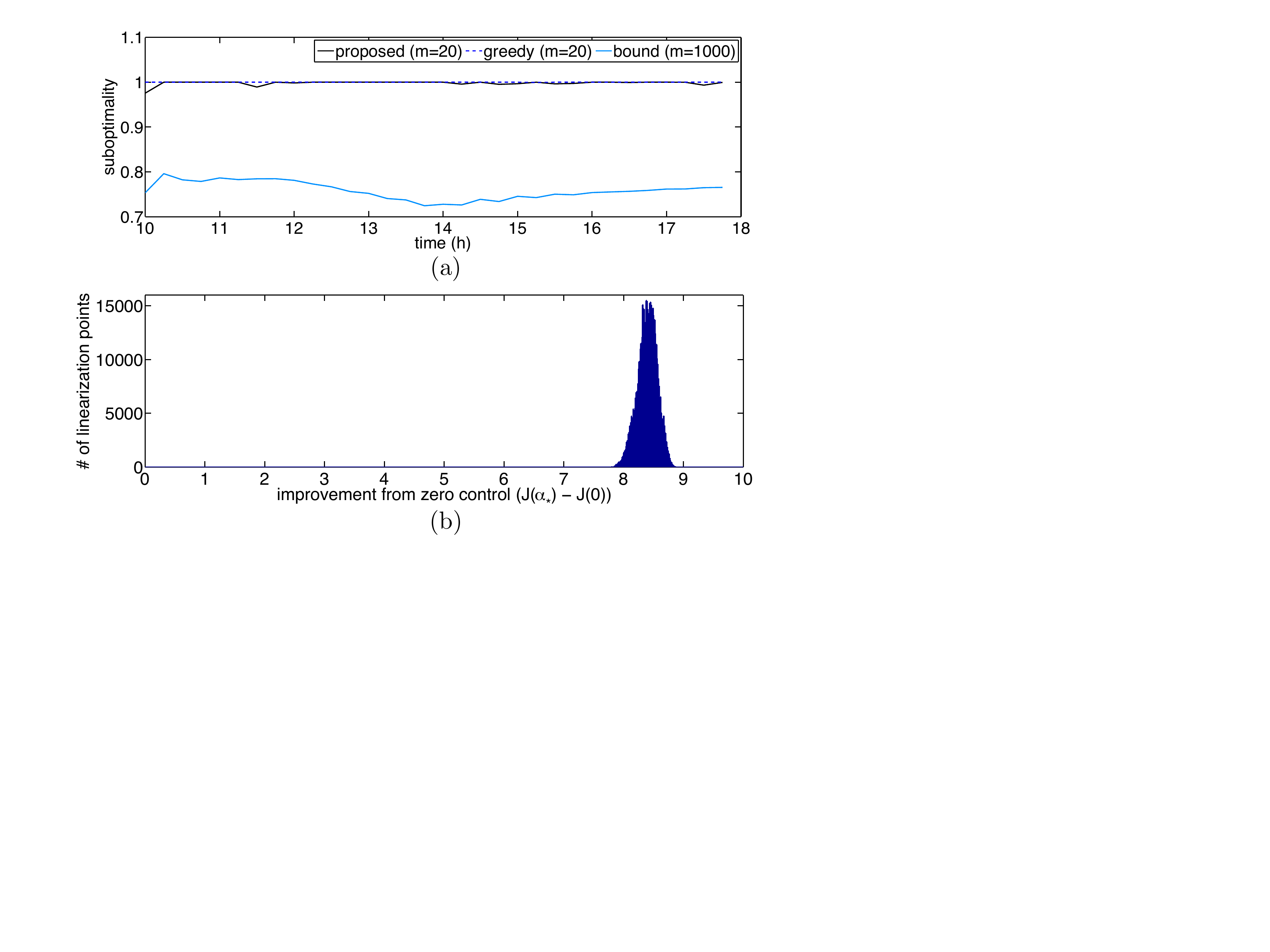}
\caption{
(a) The suboptimality bound $\rho_* = \hat{\rho}_*$ in the simulation with $m=1000$; and the performance comparison of the proposed algorithm and the greedy algorithm to the oracle when $m=20$;
 (b) Robustness test for the performance of the proposed algorithm with respect to the linearization point $\bar{\alpha}_1$ (with $m=20$).}
 \label{fig:opt_Tracking}
 \end{center}
\end{figure}

We first set the number of evaporator units to $m =1000$ and every $10$ units have the configuration in Fig.~\ref{fig:refrigerator}.
This problem approximately takes into account 25 supermarket stores.
The model parameters of the first 10 units are selected as the nominal parameter set. The model parameters of the remaining 990 units are chosen by perturbing the nominal parameter set by $\pm 10\%$ with a uniform random distribution. 
The power consumption by unit $i$ is set as $c_i = 10$kW. 
We solve the approximate problem of \eqref{opt_ex1} for $k = 1, \cdots, 32$ using the proposed algorithm.
The first five entries of the approximate solution are shown in Fig.~\ref{fig:resultTracking} (a), in which we linearize the objective function at $\bar{\alpha}^k = 0$.
The alternating pattern of the control induces that the room temperatures do not deviate significantly from $[0^\circ\mbox{C}, 4^\circ\mbox{C}]$ as shown in Fig.~\ref{fig:resultTracking} (b).
The suboptimality bound, $\rho_*$, provided in Corollary \ref{bound_improved} is computed at $k =1,\cdots, 32$. The computed values suggest that the approximate solution is at least 0.7-optimal solution for all time as shown in Fig.~\ref{fig:opt_Tracking} (a).
This suboptimality bound is better than that of the multi-linear relaxation-based local search algorithm in~\cite{Vondrak2013} for non-monotone submodular maximization with knapsack constraints, which gives at least a $(3-\sqrt{5})/2 \approx 0.309$-optimal solution.

To compute the actual suboptimality, we compare the approximate solution with the optimal solution by considering a problem with $20$ refrigerator evaporator units.
As shown in Fig.~\ref{fig:opt_Tracking} (a), the performance of the proposed approximation algorithm is greater than $95\%$ of the oracle's performance. 
In this case, the greedy algorithm performs optimally; however, we will see in the next subsection that it can get stuck at a local optimum in the presence of a more complicated constraint. 
The proposed algorithm takes 0.015 seconds to solve this problem while the greedy algorithm and exhaustive search take 0.57 seconds and 3112 seconds, respectively.

Considering $10$ refrigerator evaporator units with a single time step, i.e., $m=10$ and $k=1$, we compare the performance of the proposed algorithm and that of the greedy algorithm with $4^{10}$ initial values such that $\bold{x}_i = 2, \cdots, 5$ for $i=1, \cdots, 10$.
The ratio, $(J(\alpha_*) - J(0)) / (J(\alpha^{\mbox{\small greedy}}) - J(0))$, is within $[0.99, 1.01]$ for over 99\% of the initial values.
Lastly, we confirm that the performance of the proposed algorithm is robust with respect to the linearization point $\bar{\alpha}^1$ as shown in Fig. \ref{fig:opt_Tracking} (b) by solving the approximate problem for $m=20$ and  $k=1$ with all possible $2^{20}$ values of $\bar{\alpha}^1$.


\subsection{Case II: Customized Operation}

\begin{figure}[tb] 
\begin{center}
\includegraphics[width =3.3in]{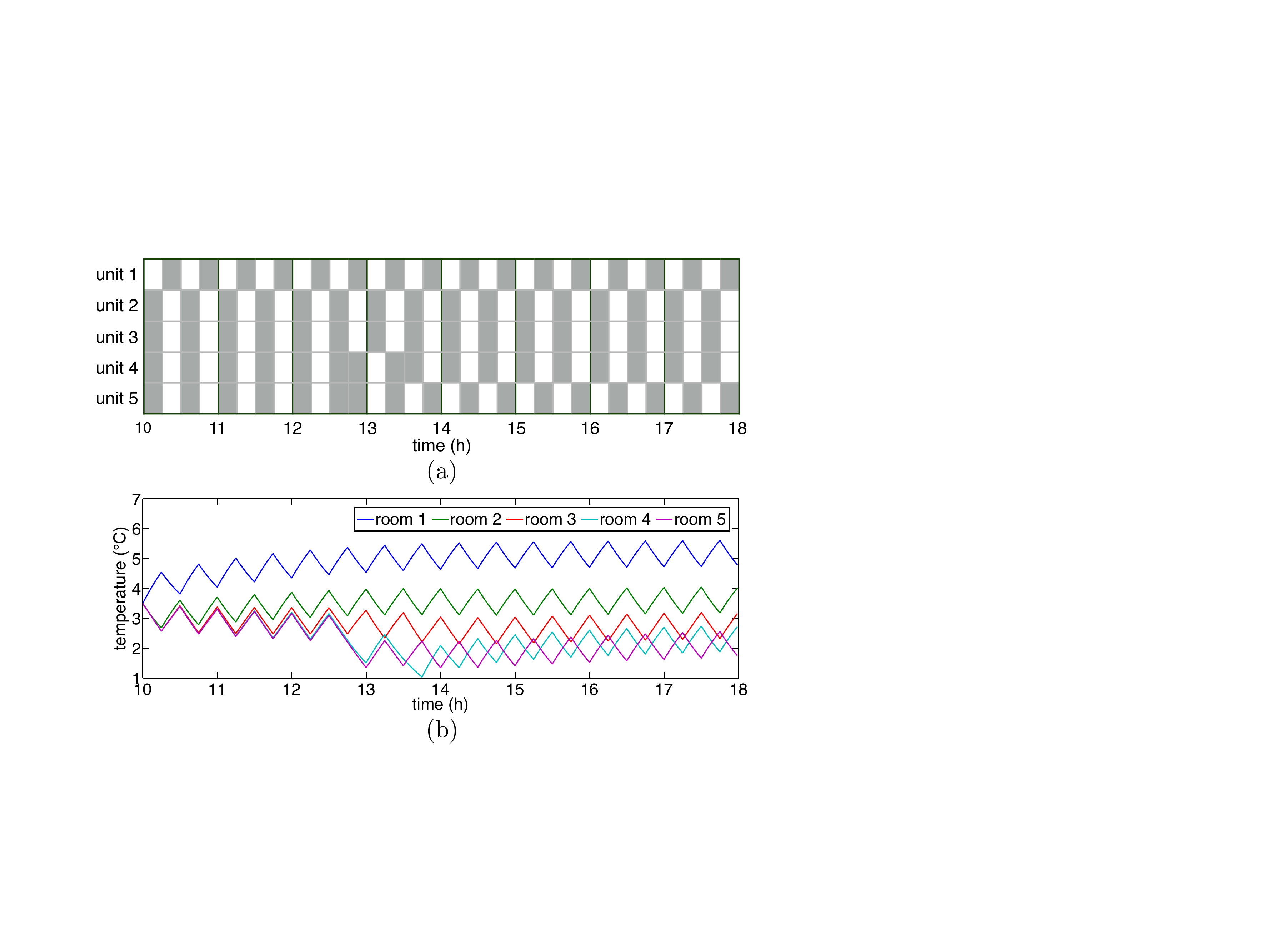}
\caption{The simulation results with approximate solution $\alpha^{k*}$, $k=1, \cdots, 32$ for
1000 units with TU constraints: (a) control signals, $\alpha^{k*}_i$, $i=1, \cdots, 5$ (grey: ON, white: OFF) (b) controlled room temperatures, $x_i$, $i=1, \cdots, 5$.}
 \label{fig:resultTU}
 \end{center}
\end{figure}

In practice, a customer may specify constraints on the operation of the refrigerators.  We consider the situation in which the constraint can be represented as 
\begin{equation} \label{TU}
\bar{Q} \alpha^k \leq \bar{r},
\end{equation}
where $\bar{Q}$ is an $m \times n$ totally unimodular (TU) matrix and $\bar{r}$ is an $m$ dimensional vector with integer entries. The usefulness of TU constraints is demonstrated in the following example.
\begin{example} \label{ex_TU}
Suppose that the power consumption by unit $1$ is comparable to the sum of the power consumptions by units $2$ and $3$. The customer has a limited budget to operate the refrigerators and therefore requests the following constraints to the aggregator:
\begin{equation}\nonumber
\begin{split}
\alpha^k_{1+10(l-1)} + \alpha^k_{2+10(l-1)}  \leq 1, &\; \; \alpha^k_{1+10(l-1)} + \alpha^k_{3+10(l-1)} \leq 1,\\
\alpha^k_{10l} + \alpha^k_{9+10(l-1)}  \leq 1, &\; \; \alpha^k_{10l} + \alpha^k_{8+10(l-1)} \leq 1
\end{split}
\end{equation}
for $l = 1, \cdots, 100$.
The rest of the units satisfy the following constraint: 
\begin{equation} \nonumber
\sum_{l=1}^{100} \sum_{j=4}^7 \alpha^k_{j+10(l-1)} \leq \bar{z}^k,
\end{equation}
where $\bar{z}_k$ is an integer.
Note that these constraints can be formulated as the inequality \eqref{TU} with a TU matrix $\bar{Q}$ and a integer vector $\bar{r}$.
\end{example}



\begin{figure}[tb] 
\begin{center}
\includegraphics[width =3.3in]{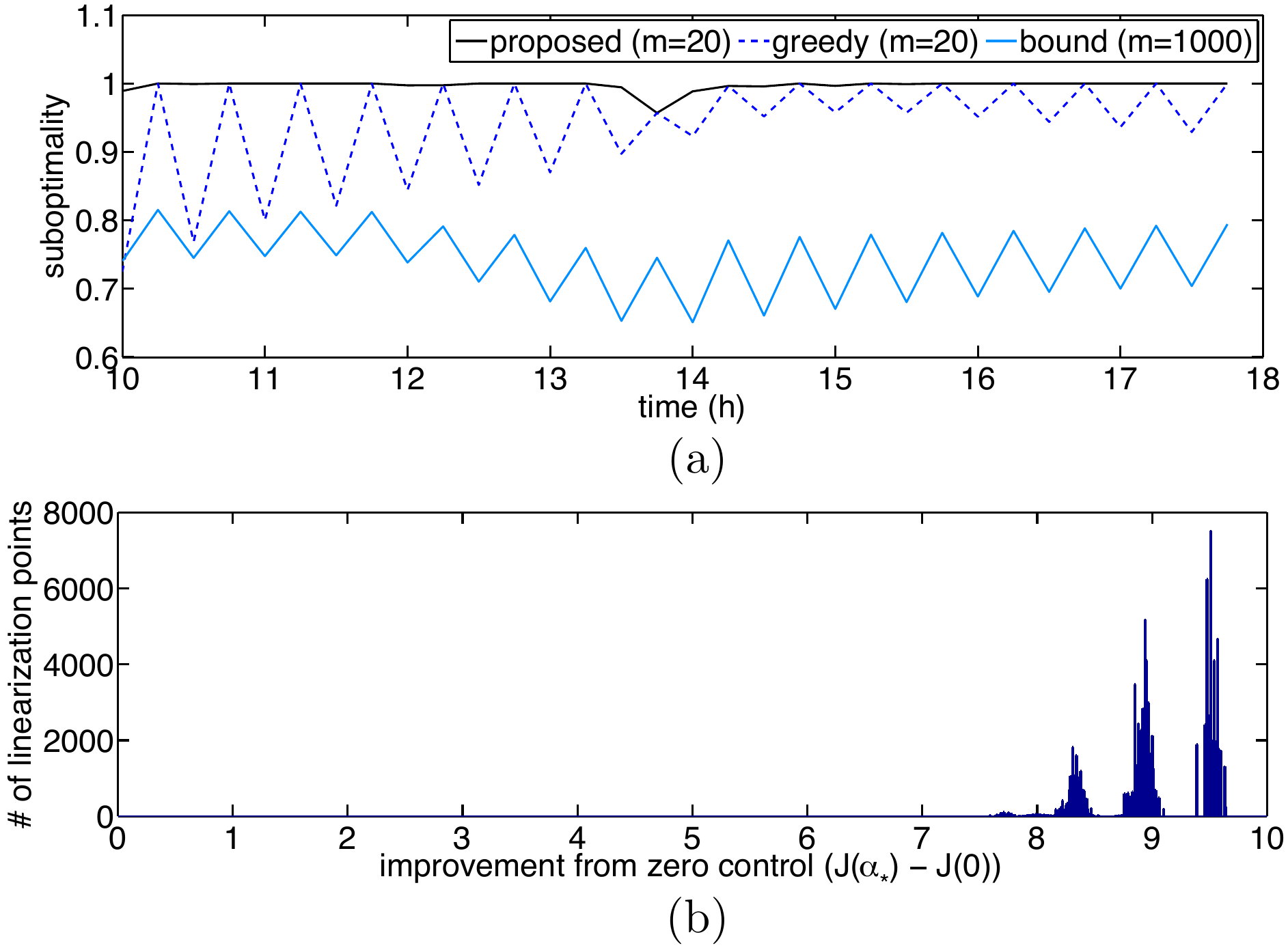}
\caption{a) The suboptimality bound $\rho_* = \hat{\rho}_*$ in the simulation with $m=1000$; and the performance comparison of the proposed algorithm and the greedy algorithm to the oracle when $m=20$;
 (b) Robustness test for the performance of the proposed algorithm with respect to the linearization point $\bar{\alpha}^1$ (with $m=20$).}
 \label{fig:opt_TU}
 \end{center}
\end{figure}

We now consider the optimization problem \eqref{opt_ex1} with TU constraint \eqref{TU} instead of \eqref{ineq_ex1}.
The decomposability of the problem depends on that of the TU constraint \eqref{TU}.
We use the same system model as that of Case I and impose the TU constraint in Example~\ref{ex_TU} with $\bar{z}^k = 2500$ for $k=9, \cdots, 16$ and $\bar{z}^k = 2000$ otherwise.
As shown in Fig.~\ref{fig:resultTU} (a), the constraints specified on units 1, 2 and 3 are satisfied by the solution of the approximate problem. 
We also note that unit 1 is used less frequently than in the previous case because if unit 1 is OFF, then units 2 and 3 can be used.
On the other hand, the greedy algorithm uses unit 1 frequently whenever the gain obtained by turning on unit 1 is greater than that by turning on unit 2 or unit 3.
This behavior is not desirable because units 2 and 3 cannot be used when  unit 1 is ON.
Therefore, the proposed algorithm outperforms the greedy algorithm as shown in Fig.~\ref{fig:opt_TU} (a).
We compute
the suboptimality bound, $\rho_*$, provided in Corollary \ref{bound_improved}  for $k =1,\cdots, 32$ as shown in Fig.~\ref{fig:opt_TU} (a).
The calculated values suggest that the approximate solution is at least 0.64-optimal solution for all time.

We again consider the problem with $20$ units to compare the approximate solution with the optimal solution. 
As shown in Fig.~\ref{fig:opt_TU}, the performance of the proposed approximation algorithm is at least 90\% of the oracle. 
On the other hand, the greedy algorithm achieves only $70- 85\%$ of the oracle 7 times out of 32.
To compare their performances with multiple initial values,
we solve the problem using the proposed approximation algorithm and the greedy algorithm for one time step, i.e., $K=1$, with $4^{10}$ initial values such that $\bold{x}_i = 2, \cdots, 5$ for $i=1, \cdots, 10$ considering $10$ refrigerator evaporator units.
The ratio, $(J(\alpha_*) - J(0)) / (J(\alpha^{\mbox{\small greedy}}) - J(0))$, is greater than $1.1$ for over 99\% of the initial values, i.e., the proposed algorithm performs at least $10\%$ better for over 99\% of the initial values. Furthermore, the average performance of the proposed algorithm is twice as high as that of the greedy algorithm.
Finally, we perform the robustness test for the proposed algorithm with respect to the liberalization point $\bar{\alpha}^1$ as in the previous subsection by solving the problem with all possible $2^{20}$ values of $\bar{\alpha}^1$. As shown in Fig. \ref{fig:opt_TU} (b), the performance does not deviate more than 15\% from its average.

\subsection{Comparison of Standard and Nonstandard Derivatives}

We now compare the performance of  approximation methods based on the standard and nonstandard derivatives. 
For comparison purpose, 
we assume that the dynamical system is given by \eqref{etp_original} for $i \in \mathcal{I}_1 := 2,4,6,8,12,14,16,18$ and
\begin{equation} \label{etp_new}
\dot{x}_i = -a_{ii} (x_i - \theta_i) - \sum_{j = 1}^n a_{ij} (x_i - x_j) - b_i  e^{-\xi_i (1 - u_i) t}
\end{equation}
for $i \in \mathcal{I}_2 := \{1,2,\cdots, 20\} \setminus \mathcal{I}_1$, where we set $\xi_i = 100$.
The modified term in the dynamical system \eqref{etp_new} models transient shutdown behavior of refrigerators after the OFF control signal is given.
For $i \in \mathcal{I}_2$, the standard and nonstandard derivatives are given by
\begin{equation} \nonumber
\begin{split}
[\derivative J(\bar{\alpha})]_i &= -b_i \xi_i \int_0^T te^{-\xi_i(1 - \bar{\alpha}_i) t} \lambda_i^{\bar{\alpha}} (t) dt,\\
[\newder J(\bar{\alpha})]_i &= -b_i \int_0^T (1 - e^{-\xi_i t}) \lambda_i^{\bar{\alpha}} (t) dt,
\end{split}
\end{equation}
respectively.
Setting $m=20$ and $k=1$,
we solve the the approximation problems based on the two derivatives with all possible $2^{20}$ values of the linearization point.
Recall that $\hat{\alpha}^*$ and $\alpha^*$ denote the solution of the approximate problems based on the nonstandard and standard derivatives, respectively.
The average of $J(\hat{\alpha}_*) - J(0)$ over all the linearization points solving the approximate problem using the nonstandard derivative
 is $13.87$, which is greater than the average $13.68$ of $J(\alpha_*) - J(0)$ obtained using the standard derivative.
Therefore,  the approximation algorithm using the nonstandard derivative performs better than that using the standard derivative on average in this problem.
This result can be explained as follows. 
Intuitively, the optimal solution should preferentially turn on refrigerator units in $\mathcal{I}_1$ because the transient behavior of refrigerator unit $i \in \mathcal{I}_2$ provides a refrigeration even when it is OFF. 
Note that the approximate solution using the nonstandard derivative preferential selects to turn on refrigerators in $\mathcal{I}_1$ because $[\hat{D}J(\bar{\alpha})]_i$  for $i \in \mathcal{I}_2$ is deflated from the case of Section \ref{target_profile}.
However, the standard derivative $[{D}J(\bar{\alpha})]_i$ is inflated for $i \in \mathcal{I}_2$ and, therefore, its approximate solution preferentially turns on refrigerators in $\mathcal{I}_2$. As a result, the approximate solution using the nonstandard derivative slightly outperforms the standard derivative. In general, the one of two approximate solutions that outperforms another is problem-dependent (see also Appendix \ref{perf_com}).

\section{Conclusion}

We have proposed approximation algorithms for optimization of combinatorial dynamical systems, in which the decision variable is a binary vector and the cost is evaluated along the solution of the systems.
The key idea of the approximation is to linearize the objective function using its derivative, which is well-defined in the feasible space of the binary decision variable. 
We proposed two different variation methods to define such derivatives. 
The approximate problem has three major advantages: $(i)$ the approximate problem is a $0$--$1$ linear program and, therefore, can be  solved by polynomial time exact or approximation algorithms;
$(ii)$ it does not require us to repeatedly solve the dynamical system;
 and $(iii)$ its solution has a provable suboptimality bound under certain concavity conditions. 
In our numerical experiments in direct load control, 
the suboptimality bound is greater than $64\%$ 
though in practice 
the performance of the proposed approximation algorithm is greater than 90\% of the oracle's performance.


\bibliographystyle{IEEEtran}

\bibliography{CDS_manuscript}

\appendices
\section{The $\epsilon$-Variational Systems}\label{evar_sys}

For given $\bar{\alpha}, {\alpha} \in \{0,1\}^m$,
the $\epsilon$-variational system associated with $(\bar{\alpha}, \alpha)$ is defined as \eqref{evar}, where its vector field is given by the convex combination of the two vector fields with $\bar{\alpha}$ and $\alpha$.
The state trajectory of \eqref{evar} is unique for given $\bar{\alpha}, {\alpha} \in \{0,1\}^m$ and is bounded on a finite time interval by Theorem 5.5 in~\cite{Bardi1997}, as shown in the following lemma.

\begin{lemma} \label{lem1}
Suppose that Assumption~\ref{a1} holds.
For any $\epsilon \in [0,1]$ and any $\bar{\alpha}, {\alpha} \in \{0,1\}^m$, the $\epsilon$-variational system \eqref{evar} associated with $(\bar{\alpha}, {\alpha})$ admits a unique solution, $x^{\epsilon (\bar{\alpha}, {\alpha})}$. In addition, $\| x^{\epsilon (\bar{\alpha}, {\alpha})} (t) \|$ is bounded by some constant independent of $\epsilon$ for all $t \in [0,T]$.
\end{lemma}

To use the $\epsilon$-variational system for defining the derivative of the payoff function, it is important to address how the $\epsilon$-variational system behaves as $\epsilon$ tends to zero  compared to the original dynamical system.
The following lemma shows that the difference $x^{\epsilon (\bar{\alpha}, {\alpha})} (t) - x^\alpha (t)$ is Lipschitz continuous in $\epsilon \in [0,1]$ for any $\bar{\alpha}, {\alpha} \in \{0,1\}^m$ (e.g., Lemma 5.6.7 in~\cite{Polak1997}).

\begin{lemma} \label{lem2}
Suppose that Assumption~\ref{a1} holds.
For any $\epsilon \in [0,1]$ and any $\bar{\alpha}, {\alpha} \in \{0,1\}^m$, there exists a constant $L$ independent of $\epsilon$ such that for all $t \in [0,T]$
\begin{equation}\nonumber
\| x^{\epsilon (\bar{\alpha}, {\alpha})} (t) - x^{\bar{\alpha}} (t) \| \leq L \epsilon.
\end{equation}
\end{lemma}

Combining the two lemmas and the dominated convergence theorem (e.g., \cite{Royden2010}), we have the following corollary.

\begin{cor} \label{cor1}
Suppose that Assumption~\ref{a1} holds.
For any $\bar{\alpha}, {\alpha} \in \{0,1\}^m$, the following equality holds:
\begin{equation}\nonumber
\lim_{\epsilon \to 0^+} \frac{1}{\epsilon} \int_0^T \| x^{\epsilon (\bar{\alpha}, {\alpha})} (t) - x^{\bar{\alpha}} (t)\|^2 dt = 0.
\end{equation}
\end{cor}

This corollary is essential to show that our proposed nonstandard derivative is well-defined and can be computed using an adjoint-based formula.

\section{Proof of Theorem~\ref{dj2}} \label{pf_dj2}

We will show a more general equality,
\begin{equation}\label{adj_form}
\begin{split}
\lim_{\epsilon  \to 0^+} &\frac{1}{\epsilon}\left [  
\mathcal{J}^{\epsilon (\bar{\alpha}, {\alpha} )} (x^{\epsilon (\bar{\alpha}, {\alpha})}) - \mathcal{J}(x^{\bar{\alpha}}, \bar{\alpha})
 \right ] = \\
 &\int_0^T \left ( 
f(x^{\bar{\alpha}}(t), \alpha) - f(x^{\bar{\alpha}}(t), \bar{\alpha}) \right )^\top \lambda^{\bar{\alpha}}(t)  + r(x^{\bar{\alpha}}(t),  \alpha) - r(x^{\bar{\alpha}}(t),  \bar{\alpha} )  dt. 
\end{split}
\end{equation}
Substituting $\alpha = \bar{\alpha} + \bold{1}_i$ and $\alpha = \bar{\alpha} - \bold{1}_i$ into the above equality, we obtain the formulae in Theorem~\ref{dj2} for $\bar{\alpha}_i=0$ and $\bar{\alpha}_i=1$, respectively.
\begin{IEEEproof}
Fix $\bar{\alpha}, {\alpha} \in\{0,1\}^m$. 
For notational simplicity, we let $\hat{x}(\cdot) := x^{\epsilon (\bar{\alpha}, {\alpha})}(\cdot) - x^{\bar{\alpha}}(\cdot)$. Then, it satisfies the following ODE:
\begin{equation}\nonumber
\begin{split}
\dot{\hat{x}}(t) &= f(\hat{x}(t) + x^{\bar{\alpha}}(t), {\bar{\alpha}}) -  f(x^{\bar{\alpha}}(t), {\bar{\alpha}}) + \epsilon ( f(\hat{x}(t) + x^{\bar{\alpha}}(t), {\alpha}) - f(\hat{x}(t) + x^{\bar{\alpha}}(t), {\bar{\alpha}}))
\end{split}
\end{equation}
with $\hat{x}(0) =0$.
The dynamical system can be rewritten as
\begin{equation} \label{perturb}
\begin{split}
\dot{\hat{x}}(t) &= \frac{\partial f(x^{\bar{\alpha}}(t), {\bar{\alpha}})}{\partial \bm{x}} \hat{x}(t) + \epsilon ( f(\hat{x}(t)+ x^{\bar{\alpha}}(t), {\alpha})  - f(\hat{x}(t) + x^{\bar{\alpha}}(t), {\bar{\alpha}})) + \sigma (\hat{x}(t), x^{\bar{\alpha}}(t)),
\end{split}
\end{equation}
where $\sigma := (\sigma_1, \cdots, \sigma_n)$ is given by
\begin{equation}\nonumber
\begin{split}
\sigma_i(\hat{\bm{x}}, \bm{x}) &:= \bold{H}_i (\hat{\bm{x}}, \bm{x}, {\bar{\alpha}}) +
 \epsilon \left ( \frac{\partial f_i(\bm{x}, {\alpha})}{\partial \bm{x}} - \frac{\partial f_i(\bm{x}, {\bar{\alpha}})}{\partial \bm{x}}  \right )\hat{\bm{x}}  + \epsilon(\bold{H}_i (\hat{\bm{x}}, \bm{x}, {\alpha}) - \bold{H}_i (\hat{\bm{x}}, \bm{x}, {\bar{\alpha}}) )
 \end{split}
\end{equation}
and $\bold{H}_i(\hat{\bm{x}}, \bm{x}, \alpha)$ denotes the higher-order terms in the Taylor expansion of $f_i(\hat{\bm{x}}+\bm{x}, \alpha)$ at $\bm{x}$, i.e., by applying the mean value theorem,
$\bold{H}_i (\hat{\bm{x}}, \bm{x}, \alpha) := \int_0^1 (1-s) \hat{\bm{x}}^\top D_{\bm{x}}^2 f_i(\bm{x} + s \hat{\bm{x}}, \alpha) \hat{\bm{x}}  ds$.
Due to Lemma~\ref{lem2} or Corollary~\ref{cor1}, for all $t \in [0,T]$,
\begin{equation}\nonumber
\lim_{\epsilon \to 0^+} \frac{1}{\epsilon} \sigma_i (\hat{x}(t), x^{\bar{\alpha}}(t))= 0.
\end{equation}

We now consider the difference
\begin{equation} \nonumber
\begin{split}
&\mathcal{J}^{\epsilon ({\bar{\alpha}}, {\alpha})} (x^{\epsilon ({\bar{\alpha}}, {\alpha})}) - \mathcal{J}(x^{\bar{\alpha}}, {\bar{\alpha}})\\
&= \int_0^T r(x^{\epsilon ({\bar{\alpha}}, {\alpha})}(t), {\bar{\alpha}}) - r(x^{\bar{\alpha}}(t), {\bar{\alpha}})  +\epsilon (r(x^{\epsilon ({\bar{\alpha}}, {\alpha})}(t), {\alpha}) - r(x^{\epsilon ({\bar{\alpha}}, {\alpha})}(t), {\bar{\alpha}}))dt + q(x^{\epsilon ({\bar{\alpha}}, {\alpha})} (T)) - q(x^{\bar{\alpha}}(T)).
\end{split}
\end{equation}
The difference can be rewritten as
\begin{equation} \nonumber
\begin{split}
&\mathcal{J}^{\epsilon ({\bar{\alpha}}, {\alpha})} (x^{\epsilon ({\bar{\alpha}}, {\alpha})}) - \mathcal{J}(x^{\bar{\alpha}}, {\bar{\alpha}})\\ 
&= \int_0^T \frac{\partial r(x^{\bar{\alpha}} (t), {\bar{\alpha}})}{\partial \bm{x}} \hat{x}(t) + \epsilon (r(x^{\bar{\alpha}}(t), {\alpha}) - r (x^{\bar{\alpha}}(t), {\bar{\alpha}})) dt + \frac{\partial q(x^{\bar{\alpha}}(T))}{\partial \bm{x}} \hat{x}(T)
+ \eta(\hat{x}, x^{\bar{\alpha}}),
\end{split}
\end{equation}
where 
\begin{equation} \nonumber
\begin{split}
\eta (\hat{x}, x):=  
&\int_0^T \bold{I}(\hat{x}(t), x(t), {\bar{\alpha}}) + \epsilon \left (
\frac{\partial r(x (t), {\alpha})}{\partial \bm{x}} -
\frac{\partial r(x (t), {\bar{\alpha}})}{\partial \bm{x}} 
\right )\hat{x}(t)\\
&\quad +\epsilon (
\bold{I} (\hat{x}(t), x(t), {\alpha}) - \bold{I} (\hat{x}(t), x(t), {\bar{\alpha}})
) dt +\bold{J} (\hat{x}(T), x(T))
\end{split}
\end{equation}
and $\bold{I}(\hat{\bm{x}}, \bm{x}, \alpha)$ and $\bold{J}(\hat{\bm{x}}, \bm{x})$  denote the higher-order terms in the Taylor expansions of $r(\hat{\bm{x}} + \bm{x}, \alpha)$ and $q(\hat{\bm{x}} + \bm{x})$ at $\bm{x}$, respectively, i.e.,
$\bold{I} (\hat{\bm{x}}, \bm{x}, \alpha) := \int_0^1 (1-s) \hat{\bm{x}}^\top D_{\bm{x}}^2 r(\bm{x} + s\hat{\bm{x}}, \alpha) \hat{\bm{x}} ds$ and
$\bold{J}(\hat{\bm{x}}, \bm{x}) := \int_0^1 (1-s) \hat{\bm{x}}^\top D_{\bm{x}}^2 q(\bm{x} + s\hat{\bm{x}}) \hat{\bm{x}} ds$.
Due to Lemma~\ref{lem2} or Corollary~\ref{cor1}, we have
\begin{equation}\nonumber
\lim_{\epsilon \to 0^+} \frac{1}{\epsilon} \eta (\hat{x}, x^{\bar{\alpha}}) = 0.
\end{equation}
Adding the inner product between the adjoint state and the system \eqref{perturb}, which is zero, to the difference, we have
\begin{equation} \label{comb1}
\begin{split}
&\mathcal{J}^{\epsilon ({\bar{\alpha}}, {\alpha})} (x^{\epsilon ({\bar{\alpha}}, {\alpha})}) - \mathcal{J}(x^{\bar{\alpha}}, {\bar{\alpha}})= \\
& \int_0^T \left (\frac{\partial r(x^{\bar{\alpha}}(t) , {\bar{\alpha}})}{\partial \bm{x}} \hat{x}(t) + \epsilon (r(x^{\bar{\alpha}}(t), {\alpha}) - r (x^{\bar{\alpha}}(t), {\bar{\alpha}})) \right ) dt + \frac{\partial q(x^\alpha (T))}{\partial \bm{x}} \hat{x}(T)\\
&+ \int_0^T (\lambda^{\bar{\alpha}} (t))^\top \left (
-\dot{\hat{x}}(t) + \frac{\partial f(x^{\bar{\alpha}}(t), {\bar{\alpha}})}{\partial \bm{x}} \hat{x}(t)  + \epsilon ( f(\hat{x}(t) + x^{\bar{\alpha}}(t), {\alpha}) - f(\hat{x}(t) + x^{\bar{\alpha}}(t), {\bar{\alpha}})) 
\right )dt\\
&+ \Theta (\hat{x}, x^{\bar{\alpha}} ),
\end{split}
\end{equation}
where 
$\Theta (\hat{x}, x^{\bar{\alpha}}) := \int_0^T (\lambda^{\bar{\alpha}})^\top \sigma (\hat{x} (t), x^{\bar{\alpha}} (t)) dt + \eta (\hat{x}, x^{\bar{\alpha}})$.
Using integration by parts, we have
\begin{equation} \label{comb2}
\begin{split}
\int_0^T (\lambda^{\bar{\alpha}})^\top \dot{\hat{x}} dt 
&= \lambda^{\bar{\alpha}}(T)^\top \hat{x} (T) -\lambda^{\bar{\alpha}}(0)^\top \hat{x} (0) 
- \int_0^T (\dot{\lambda}^{\bar{\alpha}}(t))^\top \hat{x}(t) dt\\
&=\frac{\partial q(x^{\bar{\alpha}}(T))}{\partial \bm{x}} \hat{x}(T) -\int_0^T (\dot{\lambda}^{\bar{\alpha}}(t))^\top \hat{x}(t) dt.
\end{split}
\end{equation}
Combining \eqref{comb1} and \eqref{comb2}, we obtain
\begin{equation}\nonumber
\begin{split}
\mathcal{J}^{\epsilon ({\bar{\alpha}}, {\alpha})} (x^{\epsilon ({\bar{\alpha}}, {\alpha})}) - \mathcal{J}(x^{\bar{\alpha}}, {\bar{\alpha}}) 
&=
\int_0^T \left (  
({\lambda}^{\bar{\alpha}}(t))^\top \frac{\partial f(x^{\bar{\alpha}}(t), {\bar{\alpha}})}{\partial \bm{x}}  +\frac{\partial r(x^{\bar{\alpha}}(t) , {\bar{\alpha}})}{\partial \bm{x}} + (\dot{\lambda}^{\bar{\alpha}}(t))^\top 
\right ) \hat{x}(t) dt\\
&+ \epsilon \int_0^T r(x^{\bar{\alpha}}(t), {\alpha}) - r(x^{\bar{\alpha}}(t), {\bar{\alpha}}) + (\lambda^{\bar{\alpha}}(t))^\top (f(x^{\bar{\alpha}}(t), {\alpha}) - f(x^{\bar{\alpha}}(t), {\bar{\alpha}})) dt\\
&+ \Theta(\hat{x}, x^{\bar{\alpha}}),
\end{split}
\end{equation}
where the first integral term on the right-hand side is equal to zero due to the definition of the adjoint system \eqref{adj}.
Since
\begin{equation}\nonumber
\lim_{\epsilon \to 0} \frac{1}{\epsilon} \Theta (\hat{x}, x^{\bar{\alpha}}) = 0,
\end{equation}
we obtain \eqref{adj_form} as desired.

The existence and the uniqueness of the state $x^{\bar{\alpha}}(t)$ and the adjoint state $\lambda^{\bar{\alpha}}(t)$ guarantee 
the existence and uniqueness of the nonstandard derivative.
Furthermore, the boundedness of $x^{\bar{\alpha}}(t)$ and $\lambda^{\bar{\alpha}}(t)$ for $t \in [0,T]$ imply that the nonstandard derivative is bounded.
\end{IEEEproof}

\section{Comparison of Standard and Nonstandard Derivatives} \label{comparison}

We first characterize a condition under which both derivatives are the same.
\begin{proposition} \label{prop_equiv}
Suppose that Assumptions~\ref{a1}, \ref{a2}, \ref{a3} and \ref{a4} hold.
If $f(\bm{x}, \: \cdot \:): \mathbb{R}^m \to \mathbb{R}^n$ and $r(\bm{x}, \: \cdot \:): \mathbb{R}^m \to \mathbb{R}$ are affine functions, then the two derivatives, $\derivative J$ and $\newder J$, are equivalent to each other.
\end{proposition}
\begin{IEEEproof}
Since $f(\bm{x}, \: \cdot \:)$ and $r(\bm{x}, \: \cdot \:)$ are differentiable and affine, we have
$\frac{\partial f(\bm{x}, \alpha)}{\partial \bm{\alpha}_i}  = f(\bm{x}, \bold{1}_i)$.
A similar inequality holds for $r$.
Comparing the adjoint-based formulae for $\derivative J$ and $\newder J$ in Proposition~\ref{prop1} and  Theorem~\ref{dj2}, respectively, with the assumption that $f(\bm{x}, \alpha)$ and $r(\bm{x}, \alpha)$ are affine in $\alpha$,
we deduce that the two derivatives are equivalent to each other.
\end{IEEEproof}
In general, $\derivative J$ and $\newder J$ are different from each other because they use different variation methods in their definitions.
We present a concrete example in which the standard derivative is   different the nonstandard derivative.
\begin{example} \label{ex2}
Suppose that $n=1$, $m>1$,
$f(\bm{x}, \alpha) =  \bm{x} + \sum_{i=1}^m e^{-\alpha_i}$,
$r(\bm{x}, \alpha) = \bm{x}$ and
the terminal payoff $q$ is set to be zero. Note that the vector field is not affine but additive in $\alpha$. 
Then,  the standard and nonstandard derivatives are given by
\begin{equation}\nonumber
\begin{split}
[\derivative J(\bar{\alpha})]_i &= -\int_0^T \lambda^{\bar{\alpha}}(t)   e^{-\bar{\alpha}_i}dt\\
[\newder J(\bar{\alpha})]_i &= \int_0^T \lambda^{\bar{\alpha}}(t)  (e^{-1} - e^{0})dt,
\end{split}
\end{equation}
respectively.
We notice that $[\derivative J(\bar{\alpha})]_i$ and $[\newder J(\bar{\alpha})]_i$ are not equal to each other. 
\end{example}


\subsection{Differentiability Issue}

Recall that the standard derivative $\derivative J$ requires the differentiability of $f(\bm{x}, \: \cdot \:)$ and $r(\bm{x}, \: \cdot \:)$, which can be restrictive in many applications. 
One can reformulate $f(\bm{x}, \: \cdot \:)$ and $r(\bm{x}, \: \cdot \:)$ as the following polynomials in $\alpha$ using the multi-linear polynomial extension:
\begin{equation} \nonumber
\begin{split}
\tilde{f}(\bm{x}, \alpha) &= \sum_{V \subseteq \Omega} f(\bm{x}, \mathbb{I}(V)) \prod_{i \in V} \alpha_i \prod_{i \in \Omega \setminus V} (1-\alpha_i),\\
\tilde{r}(\bm{x}, \alpha) &= \sum_{V \subseteq \Omega} r(\bm{x}, \mathbb{I}(V)) \prod_{i \in V} \alpha_i \prod_{i \in \Omega \setminus V} (1-\alpha_i),
\end{split}
\end{equation}
where $\Omega := \{1, \cdots, m\}$ and
$\mathbb{I}: 2^\Omega \to \{0,1\}^m$ is the set indicator function.
However, each of these representations requires $2^m$ calculations in the worst case.
Therefore, it is not computationally tractable to construct the multi-linear polynomial representations of $f$ and $r$.

The nonstandard derivative $\newder J$ is a good alternative to resolve this differentiability issue.
Note that this nonstandard derivative fully takes  advantage of the fact that the problem is associated with a dynamical system:
the construction of the nonstandard derivative is possible because we are able to utilize the vector field of the dynamical system as a relaxation tool.
This convex combination approach for vector fields and running payoffs naturally resolves the differentiability issue. 

\subsection{Performance Comparison}\label{perf_com}

As suggested in Section \ref{solution}, we solve the approximate problem \eqref{app_opt} twice: once  using the standard derivative $\derivative J$ and again using the nonstandard derivative $\newder J$.
Between the two approximate solutions, the solution that gives a larger payoff is chosen.
Despite this practical advantage of using the two derivative concepts, 
it is still valuable to have an insight on the comparison of their effects on the proposed approximation.
We consider a simple example, where $n = 1$, $m=2$, and the vector field and the running payoff are given by  
\begin{equation} \nonumber
\begin{split}
f(\bm{x}, \alpha) = \bm{x} + \alpha_1^3 + 2 \alpha_2, \quad r (\bm{x}, \alpha) = \bm{x}^2
\end{split}
\end{equation}
and the terminal payoff is set to be zero.
The solutions of the primal and adjoint systems are given by
$x^{\alpha}(t) = e^t\bold{x} + \int_0^t e^{t-\tau} (\alpha_1^3 + 2 \alpha_2) d\tau$ and
$\lambda^{\alpha}(t) = \int_0^{T-t} e^{T-t-\tau} 2x^\alpha (T-t-\tau) d\tau$, 
respectively.
Suppose that the initial value $\bold{x}$ is positive. Then, $x^\alpha (t)$ is positive for any $t \in [0,T]$ and for any $\alpha \in \{0,1\}^2$  and, therefore, so $\lambda^\alpha (t)$ is.
The adjoint-based formulae in Proposition~\ref{prop1} and Theorem~\ref{dj2} for the two  derivatives imply that
\begin{equation} \nonumber
\begin{split}
\derivative J (\bar{\alpha}) &= \int_0^T \begin{bmatrix} 3\bar{\alpha}_1 \\ 2 \end{bmatrix} \lambda^{\bar{\alpha}} (t) dt,\\
\newder J (\bar{\alpha}) &= \int_0^T \begin{bmatrix} 1 \\ 2 \end{bmatrix} \lambda^{\bar{\alpha}} (t) dt.
\end{split}
\end{equation}
Suppose that the constraint $\| \alpha \|_0 \leq 1$ is imposed. In this case, the optimal solution is $(0,1)$.
If we linearize the optimization problem at $\bar{\alpha} = (1,1)$, then the approximate solution based on $DJ(\bar{\alpha})$ is $(1,0)$, while that based on $\hat{D} J (\bar{\alpha})$ is $(0,1)$, which corresponds to the optimal solution.
The reason why the first derivative gives a wrong solution is that it introduces a `bias' in its first entry due to the cubic term in $\alpha_1$.
Here, we do not overstate that the second derivative performs better than the first because this bias might help  find an optimal solution in other cases.
Nevertheless, it is worth noting that the second derivative does not introduce this bias. 
We believe that this observation can stimulate further research on the performance comparison of the two derivatives and theoretic investigation on this bias in the future.

\end{document}